\documentclass[a4paper,twoside,10pt]{article} 
\usepackage{amsfonts,amssymb,mathrsfs,amsmath}
\usepackage[dvips]{color}
\usepackage{graphicx}
\usepackage{multirow}
\usepackage{tabularx}
\title{On Keller's conjecture in dimension seven}
\date{}
\author{Andrzej P. Kisielewicz and Magdalena \L ysakowska\\
\\
{\small Wydzia{\l} Matematyki, Informatyki i Ekonometrii, Uniwersytet Zielonog\'orski}\\
{\small ul. Z. Szafrana 4a, 65-516 Zielona G\'ora, Poland}\\
{\small A.Kisielewicz@wmie.uz.zgora.pl}\\
{\small M.Lysakowska@wmie.uz.zgora.pl}\\
}
\numberwithin{equation}{section}
\newtheorem{pr}{\sc Proposition}
\newtheorem{ex}[pr]{\sc Example}
\newtheorem{lemat}[pr]{\sc Lemma}
\newtheorem{st}[pr]{\sc Statement}
\newtheorem{tw}[pr]{\sc Theorem}
\newtheorem{wn}[pr]{\sc Corollary}
\newtheorem{df}{\sc Definition}

\newtheorem{uw}{\sc Remark}
\newtheorem{uwi}[uw]{\sc Remarks}

\newtheorem{nap}{\sc Example }

\newtheorem{nps}[nap]{\sc Examples}

\def\ka #1{\mathscr{#1}}
\def\kal #1 #2{\mathscr{#1}^{#2}}
\def\proof{\noindent \textit{Proof.\,\,\,}}

\def\zet{\mathbb{Z}}

\def\er{\mathbb{R}}

\def\te{\mathbb{T}}

\def\Pud #1{\operatorname{Box}(#1)}

\def\iver #1{\mbox{\tt [} #1 \mbox{\tt]}}

\begin{document}

\numberwithin{pr}{section}
\numberwithin{uw}{section}
\maketitle
\begin{abstract}
A cube tiling of $\er^d$ is a  family of pairwise disjoint cubes $[0,1)^d+T=\{[0,1)^d+t:t\in T\}$ such that $\bigcup_{t\in T}([0,1)^d+t)=\er^d$. Two cubes $[0,1)^d+t$, $[0,1)^d+s$ are called a twin pair if $|t_j-s_j|=1$ for some $j\in [d]=\{1,\ldots, d\}$ and $t_i=s_i$ for every $i\in [d]\setminus \{j\}$. 
In $1930$, Keller conjectured that in every cube tiling of $\er^d$ there is a twin pair. Keller's conjecture is true for dimensions $d\leq 6$ and false for all dimensions $d\geq 8$. For $d=7$ the conjecture is still open. Let $x\in \er^d$, $i\in [d]$, and let $L(T,x,i)$ be the set of all $i$th coordinates $t_i$ of vectors $t\in T$ such that $([0,1)^d+t)\cap ([0,1]^d+x)\neq \emptyset$ and $t_i\leq x_i$. It is known that if $|L(T,x,i)|\leq 2$ for some $x\in \er^7$ and every $i\in [7]$ or  $|L(T,x,i)|\geq 6$ for some $x\in \er^7$ and $i\in [7]$, then Keller's conjecture is true for $d=7$. In the present paper we show that it is also true for $d=7$ if $|L(T,x,i)|=5$ for some $x\in \er^7$ and $i\in [7]$. Thus, if there is a counterexample to Keller's conjecture in dimension seven, then $|L(T,x,i)|\in \{3,4\}$ for some $x\in \er^7$ and $i\in [7]$. 


\end{abstract}
\section{Introduction}

A {\it cube tiling} of $\er^d$ is a  family of pairwise disjoint cubes $[0,1)^d+T=\{[0,1)^d+t:t\in T\}$ such that $\bigcup_{t\in T}([0,1)^d+t)=\er^d$. Two cubes $[0,1)^d+t$, $[0,1)^d+s$ are called a {\it twin pair} if $|t_j-s_j|=1$
for some $j\in [d]=\{1,\ldots, d\}$ and $t_i=s_i$ for every $i\in [d]\setminus \{j\}$. In $1907$, Minkowski \cite{Min} conjectured that in every {\it lattice} cube tiling of $\er^d$, i.e. when $T$ is a lattice in $\er^d$, there is a twin pair, and in $1930$, Keller \cite{Ke1} generalized this conjecture to arbitrary cube tiling of $\er^d$. Minkowski's conjecture was confirmed by Haj\'os \cite{H} in $1941$. In $1940$, Perron \cite{P} proved that Keller's conjecture is true for all dimensions $d\leq 6$. 
 
In 1992, Lagarias and Shor \cite{LS1}, using ideas from  Corr\'adi's and Szab\'o's papers \cite{CS2,Sz2}, constructed a cube tiling of $\er^{10}$ which does not contain a twin pair  and thereby refuted  Keller's cube tiling conjecture. In $2002$, Mackey \cite{M} gave a counterexample to Keller's conjecture in dimension eight, which also shows that this conjecture is false in dimension nine. For $d=7$ Keller's conjecture is still open.  

Let  $[0,1)^d+T$ be a cube tiling, $x\in \er^d$ and $i\in [d]$, and let $L(T,x,i)$ be the set of all $i$th coordinates $t_i$ of vectors $t\in T$ such that $([0,1)^d+t)\cap ([0,1]^d+x)\neq \emptyset$ and $t_i\leq x_i$ (Figure 1). 
For every cube tiling $[0,1)^d+T$, $x\in \er^d$ and $i\in [d]$ the set $L(T,x,i)$ contains at most $2^{d-1}$ elements.

\bigskip
{\center
\includegraphics[width=5cm]{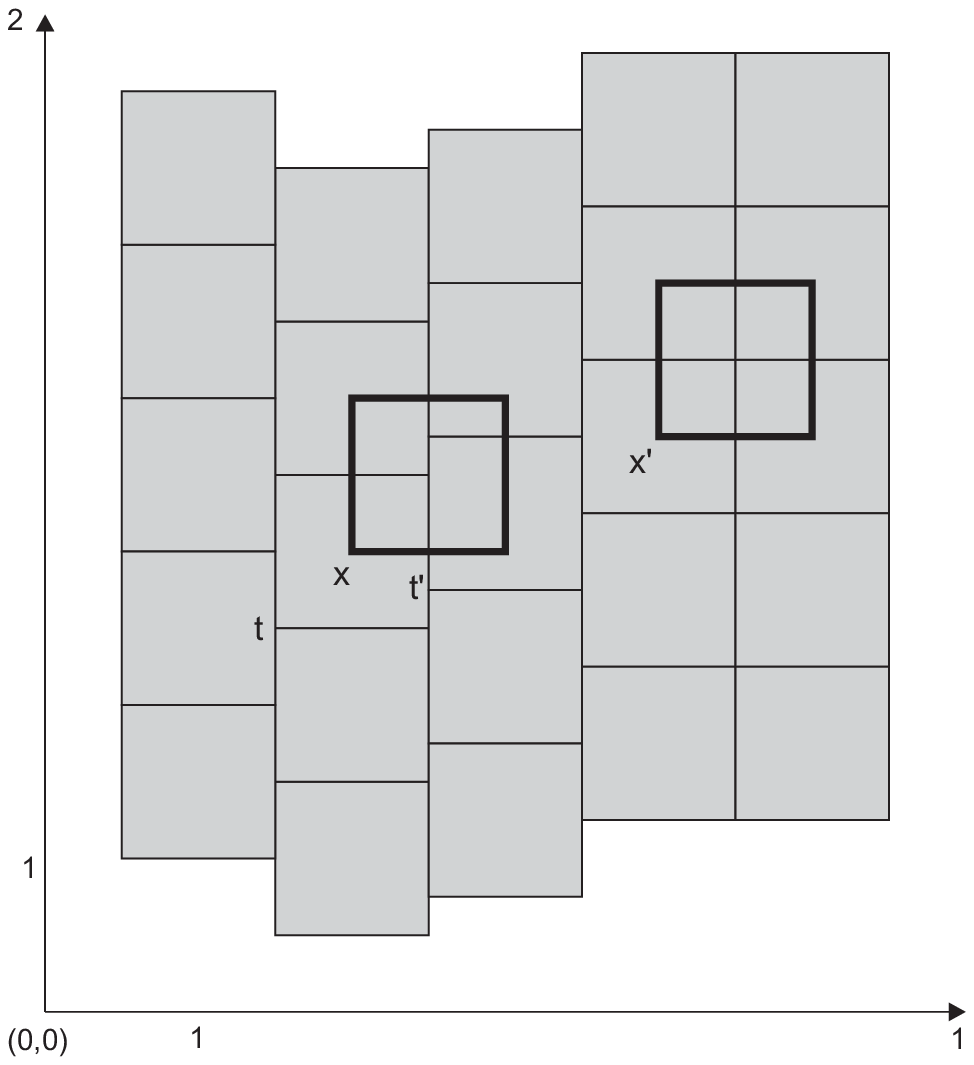}\\
}

\noindent{\footnotesize 
Fig. 1. A portion of a cube tiling $[0,1)^2+T$ of $\er^2$. The number of elements in $L(T,x,i)$ depends on the position of $x\in \er^2$. For $x=(2,3)$, we have $L(T,x,1)=\{3/2\}(=\{t_1\})$ and $L(T,x,2)=\{5/2,11/4\}(=\{t_2,t_2'\})$, while for $x'=(4,15/4)$, we have $L(T,x',1)=\{7/2\}$ and $L(T,x',2)=\{13/4\}$.
  }

\medskip

In 2010, Debroni et al. \cite{De} computed, using the supercomputer Cray XT5 Kraken, that Keller's conjecture is true for all cube tilings $[0,1)^7+T$ of $\er^7$ such that  $T\subset (1/2)\zet^7$. This result shows that Keller's conjecture is true for cube tilings of $\er^7$ with $|L(T,x,i)|\leq 2$ for some $x\in \er^7$ and every $i\in [7]$ (see \cite[Section 1]{Kis}). 
In \cite{Kis} we showed that Keller's conjecture is true for cube tilings $[0,1)^7+T$ of $\er^7$ for which $|L(T(x,i))|\geq 6$ for some $x\in \er^7$ and  $i\in [7]$.
In this paper we prove that 

\medskip
{\it Keller's conjecture is true for cube tilings $[0,1)^7+T$ of $\er^7$ for which 

 $|L(T(x,i))|=5$ for some $x\in \er^7$ and  $i\in [7]$.}

\medskip
It follows from the above results that if there is a counterexample to Keller's conjecture in dimension seven, then $|L(T,x,i)|\in \{3,4\}$ for some $x\in \er^7$ and $i\in [7]$. 
  
Keller's cube tiling conjecture was not as clearly motivation as Minkowski's conjecture was. Recall that, the existence of a twin pair in a lattice tiling $[0,1)^d+T$ determines the form of a basis for the lattice $T$. Keller's conjecture was rather a generalization of Minkowski's conjecture. In our opinion the paper of Lagarias and Shor \cite{LS2} presents the problem of the existence of twin pairs in cube tilings of $\er^d$ in an appropriate manner. In this excellent work, roughly speaking, the authors measure distances between some of the cubes in a tiling $[0,1)^d+T$. 
When the dimension of the space increases, the distances between cubes can also increase. In \cite{LS2} Lagarias and Shor gave  an estimation of how fast the distances between cubes increase. A twin pair is a pair of cubes with the minimal possible distance in a cube tiling. It follows from Perron's result that for $d\leq 6$ in an arbitrary cube tiling of $\er^d$ there are cubes which are closed (twin pairs). From Mackey's example we know that in dimension eight the process of cubes moving away in cube tilings has started. Resolving Keller's conjecture for $d=7$ will answer the question whether this process had already begun in dimension seven.    
  
Working on Keller's conjecture has, on the one hand, provided the opportunity of answering an old query in tiling theory, and on the other hand is the beginning of a much deeper and more interesting investigation into the structure of cube tilings of $\er^d$ in the spirit of Lagarias's and Shor's ideas contained the paper  \cite{LS2}. These investigations, besides describing the structure of tilings, can provide new tools which can be used in various areas of combinatorics. For example, in \cite{Kis1} we showed how a cube tiling code designed in \cite{LS2} can be used to obtain an interesting partitions and matchings of a $d$-dimensional cube. Moreover, the result of Debroni et al.\cite{De} and \cite[Theorem 3.1]{Kis} give a new proof of Keller's conjecture in dimensions $d\leq 6$ (see Remark \ref{Kel}).  
  

The presented paper is lengthy, but much of the content is in the form of a summary of results in the form of tables and figures. The proof of the crucial result (Theorem \ref{s12}) that allows us to prove the assertion on Keller's conjecture for $|L(T(x,i))|=5$  is based on computations, and these need reductions. The two longest and most arduous sections of the  paper, Section 3 and 4, contain the preparatory results for those subsequent reductions, and most of  Section 5 consists of the presentation of the initial data for  these computations. Hence, the reader who would first like an overview of how the discussed case of Keller's conjecture is proven may skip these sections and continue reading from Theorem \ref{s12}. 

To make the paper self-contained we have collected the basic notions in Section 2. We use a very abstract language, in the form of systems of abstract words, but in the long run such an approach simplifies the reasoning. Moreover, there is a nice "model" of systems of words: this is a family of pairwise disjoint translates of the unit cube in the flat torus  $\te^d=\{(x_1,\ldots ,x_d)({\rm mod} 2):(x_1,\ldots ,x_d)\in \er^d\}$. We will present the content of the paper on the example of the tilings of $\te^d$ by unit cubes. 

A set $F\subset \te^d$ is called a {\it polycube} if $F$ has a { \it tiling} by translates of the unit cube, i.e., there is a family of pairwise disjoint translates of the unit cube $[0,1)^d+T$, $T\subset \te^d$, such that $\bigcup_{t\in T}[0,1)^d+t=F$. Clearly, if $[0,1)^d+T$ and $[0,1)^d+T'$ are tilings of $F$, then $|T|=|T'|$. The question of how many tilings the polycube F has is a basic in tiling theory. As we show in Section 2 the case $|L(T(x,i))|=5$ is reduced to the following classic problem: {\it For a polycube $F$ which has a tiling consisting of 12 cubes such that no two cube form a twin pair, find all tilings of $F$ by translates of the unit cubes assuming that $F$ has at least two disjoint tilings without twin pairs.} The case  $|L(T(x,i))|\geq 6$ resolved in \cite{Kis} relies on showing that no  two such tilings exist for a polycube with 11 cubes or less. 

In  graph theory knowing the structure of  small graphs (graphs with a few vertices) plays an important role. Similarly in cube tilings, it is absolutely necessary to know the structures of all tilings of polycubes $F$ with a few cubes. These small polycubes are described in Section 3.   
In Section 4 we establish necessary conditions that have to be fulfilled by the above mentioned tilings with 12 cubes.

In Section 5, based on the results from the previous two sections, we first collect all the initial configurations of cubes for the computations. These configurations are necessary for us to make the computations as the number of all cases that would have to be considered by the computer program is more than $\binom{64}{12}3^{72}$.  At the end of Section 5 based on the results of the computations we  give  proof of the theorem on the structure of tilings of the polycube $F$ with 12 cubes without twin pairs (Theorem \ref{s12}).  Finally, in Section 6 using Theorem \ref{s12} we prove in an easy manner that Keller's cube tiling conjecture is true for tilings  $[0,1)^7+T$ of $\er^7$ with $|L(T(x,i))|=5$ for some $x\in \er^7$ and  $i\in [7]$.

\section{Basic notions}

In this section we present the basic notions on dichotomous boxes and words (details can be found in \cite{GKP,KP}). We start with systems of boxes.

A non-empty set $K \subseteq X=X_1\times\cdots \times X_d$ is called a \textit{ box} if $K=K_1\times\cdots \times K_d$ and
$K_i\subseteq X_i$ for each $i\in [d]$. By  $\Pud {X}$ we denote the set of all boxes in $X$. 
The set $X$ will be called a $d$-{\it box}. 
The box  $K$ is said to be \textit{ proper} if $K_i\neq X_i$ for each $i\in [d]$.  
Two boxes $K$ and $G$ in $X$ are called \textit{dichotomous}
if there is $i\in [d]$ such that $K_i=X_i\setminus G_i$. A \textit{suit} is any collection of pairwise
dichotomous boxes. A suit is \textit{proper} if it consists of proper boxes. 
A non-empty set $F\subseteq X$ is said to be a \textit{ polybox} if
there is a suit $\ka F$ for $F$, i.e. if $\bigcup \ka F=F$. A polybox $F$ is {\it rigid} if it has exactly one suit. (Figure 7 presents the suit for a rigid polybox. The polyboxes  $\bigcup \ka F^{3,A}$ and $\bigcup \ka F^{3,A'}$ in Figure 2 are not rigid).




The important property of proper suits is that, for every proper suits $\ka F$ and $\ka G$ for a polybox $F$, we have $|\ka F|=|\ka G|$ (this property is obvious for two tilings of a polycube $F\subset \te^d$ but not for polyboxes). Thus, we can define a {\it box number} $|F|_0=$ the number of boxes in any proper suit for $F$ (in Figure 2 we have $|\bigcup \ka F^{3,A}|_0=3$). A proper suit for a $d$-box $X$ is called a {\it minimal partition } of $X$ (Figure 2). Every minimal partition of a $d$-box has $2^d$ boxes. 


Two boxes $K,G\subset X$ are said to be a {\it twin pair} if $K_j=X_j\setminus G_j$ for some $j\in [d]$ and $K_i=G_i$ for every $i\in [d]\setminus\{j\}$. Observe that, the suit for a rigid polybox can not contain a twin pair.

 
Every two cubes $[0,1)^d+t$ and $[0,1)^d+p$ in an arbitrary cube tiling $[0,1)^d + T$ of $\er^d$ satisfy {\it Keller's condition} (\cite{Ke1}): there is $i\in [d]$ such that $t_i-p_i\in \zet\setminus\{0\}$, where $t_i$ and $p_i$ are $i$th coordinates of the vectors $t$ and $p$. For a cube $[0,1]^d+x$, where $x=(x_1,...,x_d)\in \er^d$, the family $\ka F_x=\{([0,1)^d+t)\cap ([0,1]^d+x)\neq\emptyset:t\in T\}$ is a partition of the cube $[0,1]^d+x$, in which, because of Keller's condition, every two boxes $K,G\in \ka F_x$ are dichotomous, i.e. there is $i\in [d]$ such that $K_i$ and $G_i$ are disjoint and $K_i\cup G_i=[0,1]+x_i$. Moreover, since cubes in cube tilings are half-open, every box $K\in \ka F_x$ is proper, and consequently the family $\ka F_x$ is a minimal partition of $[0,1]^d+x$. The structure of the partition $\ka F_x$ reflects the local structure of the cube tiling $[0,1)^d+T$. Obviously, a cube tiling $[0,1)^d+T$ contains a twin pair if and only if the partition $\ka F_x$ contains a twin pair for some $x\in \er^d$ (see Figure 1). 

\subsection{Our approach}

Below we sketch our approach to the problem of the existence of twin pairs in a cube tiling $[0,1)^7+T$ of $\er^7$ with $|L(T,x,i)|=5$. To do this we describe the structure of a minimal partition.
A graph-theoretic description of this structure can be found in \cite{CS1} (see also \cite{LS2}).  
  
Let $X$ be a $d$-box. A set $F\subseteq X$ is called an {\it $i$-cylinder} if 
$$
l_i\cap F=l_i \; \; \; {\rm or}\; \; \; l_i\cap F=\emptyset,
$$
where $l_i=\{x_1\}\times \cdots \times \{x_{i-1}\}\times X_i \times \{x_{i+1}\}\times \cdots \times \{x_d\}$ and $x_j\in X_j$ for $j\in [d]\setminus\{i\}$ (Figure 2). 

Let $\ka F$ be a minimal partition, and 
let $A\subset X_i$ be a set such that there is a box $K\in \ka F$ with $K_i\in \{A, X_i\setminus A\}$.  Let
$$
\ka F^{i,A}=\{K\in \ka F: K_i=A\} \; \; {\rm and} \; \; \ka F^{i,A'}=\{K\in \ka F: K_i=X_i\setminus A\}.
$$
Since the boxes in $\ka F$ are pairwise dichotomous, the set $\bigcup (\ka F^{i,A}\cup \ka F^{i,A'})$ is an $i$-cylinder, and the set of boxes $\ka F^{i,A}\cup \ka F^{i,A'}$ is a suit for it. As $|\ka F|=2^d$, it follows that the boxes in $\ka F$ can form at most $2^{d-1}$ pairwise disjoint $i$-cylinders. More precisely, for every $i\in [d]$ there are sets $A^1,\ldots ,A^k\subset X_i$ such that $A^n\not\in \{A^m,X_i\setminus A^m\}$ for every $n,m\in [k], n\neq m$, and 
$$
\ka F=\ka F^{i,A^1}\cup \ka F^{i,(A^1)'}\cup \cdots \cup \ka F^{i,A^k}\cup \ka F^{i,(A^k)'}. 
$$
The boxes in $\ka F$ are proper, and hence $|\ka F^{i,A^n}\cup \ka F^{i,(A^n)'}|\geq 2$ for every $i\in [d]$ and $n\in [k]$. Thus $k\leq 2^{d-1}$ and consequently $|L(T,x,i)|\leq 2^{d-1}$ for every cube tiling $[0,1)^d+T$, $x\in \er^d$ and $i\in [d]$, as $|L(T,x,i)|$ is the number of all $i$-cylinders in the partition $\ka F_x$.  

If $K$ is a box in $X$, $\ka G$ is a family of boxes, $x\in X$ and $i\in [d]$, then  
$$
(K)_{i}=K_1\times\cdots \times K_{i-1}\times K_{i+1}\times\cdots \times K_d,
$$
$$
(\ka G)_{i}=\{(K)_{i}:K\in \ka G\}
$$
and 
$$
(x)_i=(x_1,\ldots ,x_{i-1},x_{i+1},\ldots ,x_d).
$$

\medskip
{\center
\includegraphics[width=10cm]{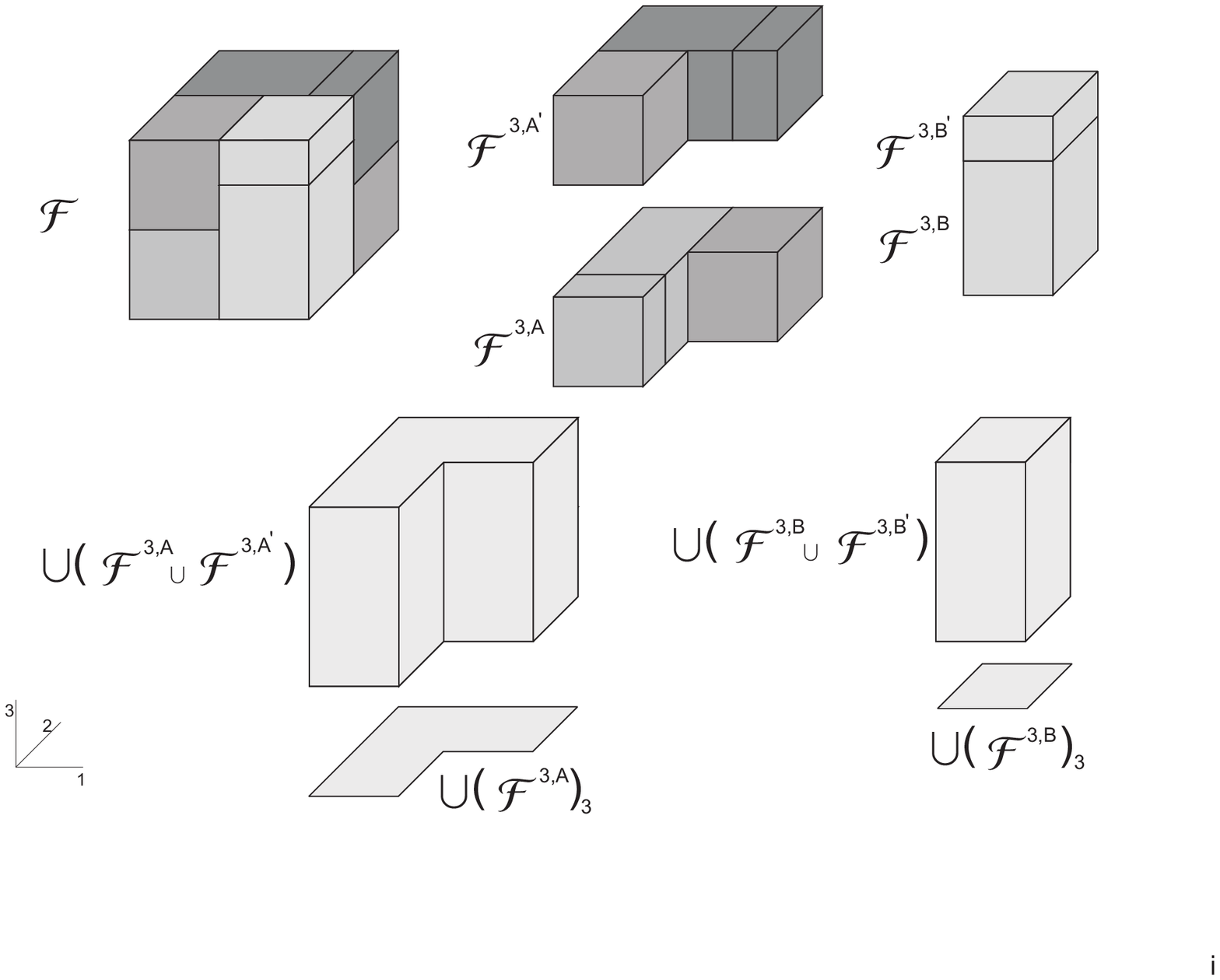}\\
}

\noindent{\footnotesize 
Fig. 2. The minimal partition $\ka F=\ka F^{3,A}\cup \ka F^{3,A'}\cup \ka F^{3,B}\cup \ka F^{3,B'}$ of the $3$-box $X=[0,1]^3$ ($A=[0,1/2)$, $B=[0,3/4)$), two $3$-cylinders and its suits.
  }
  
\medskip
Since $\bigcup (\ka F^{i,A}\cup \ka F^{i,A'})$ is an $i$-cylinder, the sets of boxes $(\ka F^{i,A})_{i}$ and $(\ka F^{i,A'})_{i}$ are two suits for the polybox $\bigcup (\ka F^{i,A})_{i}=\bigcup (\ka F^{i,A'})_{i}$, which is a polybox in the $(d-1)$-box $(X)_{i}$ (Figure 2). Note that, as $(\ka F^{i,A})_{i}$ and $(\ka F^{i,A'})_{i}$ are proper suits for the polybox $\bigcup (\ka F^{i,A})_{i}$, we have  $|(\ka F^{i,A})_{i}|=|(\ka F^{i,A'})_{i}|$ 

Let now $[0,1)^7+T$ be a cube tiling of $\er^7$, and let $\ka F_x$ be as defined above. If $|L(T,x,i)|=5$ for some $i\in [7]$, then
$$
\ka F_x=\ka F^{i,A^1}\cup \ka F^{i,(A^1)'}\cup \ldots \cup \ka F^{i,A^5}\cup \ka F^{i,(A^5)'}. 
$$
Assume that there are no twin pairs in the tiling $[0,1)^7+T$. Then $\ka F_x$ does not contain a twin pair. It follows from  \cite[Theorem 3.1]{Kis} (see Theorem \ref{12} in Section 3) that $|\ka F^{i,A^k}|\geq 12$ for every $k\in [5]$. Thus, there is at least one $k\in [5]$ such that $|\ka F^{i,A^k}|=12$ because $|\ka F_x|=128$ and $|\ka F^{i,A^k}|=|\ka F^{i,(A^k)'}|$ for every $k\in [5]$. The main effort in the paper will be rely on describing the structure of all twin pairs free suits $\ka F^{i,A^k}\cup \ka F^{i,(A^k)'}$ for $i$-cylinders  $\bigcup (\ka F^{i,A^k}\cup \ka F^{i,(A^k)'})$ such that $|\ka F^{i,A^k}|=12$. Knowing this structure, we will be able to prove that Keller's conjecture is true for a cube tiling  $[0,1)^7+T$ with $|L(T,x,i)|=5$. 

\begin{uw} 
\label{Kel}
{\rm If $d=6$ and $\ka F_x$ has at least three $i$-cylinders for some $i\in [6]$, then  $|\ka F^{i,A^k}|\leq 10$ for some $k$, and by \cite[Theorem 3.1]{Kis}, there is a twin pair in $\ka F^{i,A^k}\cup \ka F^{i,(A^k)'}$. If $\ka F_x$ has at most two $i$-cylindrs for every $i\in [6]$, then it follows from the result of Debroni et al.\cite{De} that there is a twin pair in $\ka F_x$. Thus, these two  results in \cite{De,Kis} give a new proof of Keller's conjecture in dimensions $d\leq 6$.  
}
\end{uw}

\subsection{Dichotomous words}

Two dichotomous boxes $K, G\subset X$ are of the forms: $K_1\times \cdots \times K_i\times \cdots \times K_d$ and $G_1\times \cdots \times (X_i\setminus K_i)\times \cdots \times G_d$ for some $i\in [d]$. To shorten this notation we can write $k_1\cdots k_i\cdots k_d$ and  $g_1\cdots k'_i\cdots g_d$. Thus, we can encode a system of dichotomous boxes as a system of {\it words} (see also \cite{LS2}). 
 Below we collect basic notions concerning dichotomous words (details can be found in \cite{KP}).

A set $S$ of arbitrary objects will be called an \textit{alphabet}, and the elements of $S$ will be called \textit{letters}. A permutation $s\mapsto s'$ of the alphabet $S$ such that $s''=(s')'=s$ and $s'\neq  s$ is said to be a \textit{complementation}. We add an extra letter $\ast$ to the set $S$ and the set $S\cup \{*\}$ is denoted by $*S$. We set $*'=*$ (the star is the only letter with this property).  Each sequence of letters $s_1\cdots s_d$ from the set $*S$ is called a \textit{word}. The set of all words of length $d$ is denoted by $(*S)^d$, and by $S^d$ we denote the set of all words $s_1\cdots s_d$ such that $s_i\neq *$ for every $i\in [d]$.  Two words $u=u_1\cdots u_d$ and $v=v_1\cdots v_d$ are \textit{dichotomous} if there is $j\in[d]$ such that  $u_j\neq *$ and $u'_j=v_j$. If $V\subset (*S)^d$ consists of pairwise dichotomous words, then we call it a \textit{polybox code} (or \textit{polybox genome}).
Two words $u,v\in (*S)^d$ are \textit{a twin pair } if there is $j\in[d]$ such that $u'_j=v_j$, where $u_j\neq *$ and  $u_i=v_i$ for every $i\in [d]\setminus\{j\}$. 

If $A\subset [d]$ and $A^c=[d]\setminus A=\{i_1<\cdots <i_n\}$, then $(u)_{A}=u_{i_1}\cdots u_{i_n}$ and  $(V)_{A}=\{(v)_{A}:v\in V\}$ for $V\subset (*S)^d$. If $A=\{i\}$, then we write $(u)_{i}$ and $(V)_{i}$ instead of $(u)_{\{i\}}$ and $(V)_{\{i\}}$, respectively. If $V\subset (*S)^d$, $s\in *S$ and $i\in [d]$, then let $V^{i,s}=\{v\in V:v_i=s\}$. If $V$ is a polybox code, then the representation
$$
V=V^{i,l_1}\cup V^{i,l_1'}\cup \ldots \cup V^{i,l_k}\cup V^{i,l_k'},  
$$
where $l_j,l_j'\in *S$ for $j\in [k]$, will be called a {\it distribution of words in} $V$.

Suppose now that for each $i\in[d]$ a mapping $f_i\colon *S\to \Pud {X_i}$ is such that $f_i(s')=X_i\setminus f_i(s)$ for $s\neq *$ and $f_i(*)=X_i$. We define the mapping $f\colon (*S)^d\to \Pud  X$ by
$$
f(s_1\cdots s_d)=f_1(s_1)\times\cdots\times f_d(s_d).
$$ 

About such defined $f$ we will say that it \textit{preserves dichotomies}. If $V\subseteq (*S)^d$, then the set of boxes $f(V)=\{f(v)\colon v\in V\}$ is said to be a \textit{realization} of the set of words $V$. Clearly, if $V$ is a polybox code, then $f(V)$ is a suit for $\bigcup f(V)$. The realization is said to be \textit{exact} if for each pair of words  
$v,w \in V$,  if  $v_i\not\in \{w_i, w'_i\}$, then $f_i(v_i)\not\in\{f_i(w_i), X_i\setminus f_i(w_i)\}$. 

A polybox code $V\subset (*S)^d$ is called  a {\it partition code } if arbitrary realization $f(V)$ of $V$ is a suit for a $d$-box $X$. Observe that, if  $V\subset S^d$ is  a partition code, then $f(V)$ is a minimal partition. 


We will exploit some abstract but very useful realization of polybox codes. This sort of realization was invented in \cite {ABHK} (but our improved construction comes from \cite{KP}), where it was the crucial tool in proving the main theorem of that paper. 
   
Let $S$ be an alphabet with a complementation, and let
$$
 ES=\{B\subset S\colon |\{s,s'\}\cap B|=1, \text{whenever $s\in S$}\},
$$
$$ 
E s=\{B\in ES\colon s\in B\}\;\; {\rm and}\;\; E*=ES.
$$   
Let $V\subset (*S)^{d}$ be a polybox code, and let $v\in V$. The {\it equicomplementary} realization of the word $v$ is the box   
$$
\breve{v}=Ev_1\times \cdots \times Ev_d
$$
in the $d$-box $(ES)^d =ES\times \cdots \times ES.$ The equicomplementary realization of the code $V$ is the family
$$
E(V)=\{\breve{v}:v\in V\}.
$$ 
If $S$ is finite, $s_1,\ldots , s_n\in S$ and $s_i\not\in\{s_j, s'_j\}$ for every $i\neq j$, then
\begin{equation}
\label{dkostki}
|E s_1 \cap\dots\cap E s_n|=(1/2^{n})|ES|.
\end{equation}
In the paper we will assume that $S$ is finite.
The value of the realization $E(V)$, where $V\subset S^d$,  lies in the above equality (which does not hold for translates of the unit interval $[0,1)$ in $\te^1$). In particular, boxes in $E(V)$ are of the same size; for $w\in E(V)$ we have $|\breve{w}|=(1/2^d)|ES|^d$. Thus, two boxes $\breve v, \breve w\subset (ES)^d$ are dichotomous if and only if $\breve v\cap \breve w=\emptyset.$ The same is true for cubes in a cube tiling of a polycube $F\subset \te^d$ and therefore working with the boxes $\breve{v}, v\in V$, we can think of them as translates of the unit cube in $\te^d$.   

Moreover, from (\ref{dkostki}) we obtain the following important lemma which was proven in \cite{Kis}: 

\begin{lemat}
\label{=c}
Let $w,u,v\in S^d$, and let $\ka D$ be a simple partition of the $d$-box $\breve w$. If boxes $\breve w\cap \breve u$ and  $\breve w\cap \breve v$ belong to $\ka D$, then there is a simple partition code $C\subset S^d$ such that $u,v\in C$. In particular, if $\breve w\cap \breve u$ and  $\breve w\cap \breve v$ form  a twin pair, then $u$ and $v$ are a twin pair.
 \hfill{$\square$}
\end{lemat}

\medskip
Let $V,W\subset (*S)^d$ be polybox codes, and let $v\in (*S)^d$. We say that $v$ is \textit{covered} by $W$, and write $v\sqsubseteq W$, if $f(v)\subseteq \bigcup f(W)$ for every mapping $f$ that preserves dichotomies. If  $v\sqsubseteq W$ for every $v\in V$, then we write $V \sqsubseteq W$.

Polybox codes $V,W\subset (*S)^d$ are said to be {\it equivalent} if $V \sqsubseteq W$ and  $W \sqsubseteq V$.  A polybox code $V\subset S^d$ is called {\it rigid} if there is no code $W\subset S^d$ which is equivalent to $V$ and $V\neq W$. Observe that, rigid polybox codes can not contain a twin pair.

\medskip
Let  
$g\colon S^d\times S^d\to \zet$ be defined by the formula
\begin{equation}
\label{g}
g(v,w)=\prod^d_{i=1}(2[v_i=w_i]+[w_i\not\in\{v_i, v'_i\}]),
\end{equation}
where $\iver {p}=1$ if the sentence $p$ is true and $\iver {p}=0$ if it is false.

Let $w\in S^d$, and let $V\subset S^d$ be a polybox code. Then 

\begin{equation}
\label{2d}
\breve{w}\subseteq \bigcup E(V)\Leftrightarrow w \sqsubseteq V \Leftrightarrow \sum_{v\in V} g(v,w)= 2^d.
\end{equation}

It follows from the definition of equivalent polybox codes $V,W\subset S^d$ and (\ref{2d}) that $V$ and $W$  are equivalent if and only if $\bigcup E(V)=\bigcup E(W)$.


\medskip
Let $s_*=*\cdots *\in (*S)^d$ and let $\bar{g}(\cdot,s_*)\colon (*S)^d\to \zet$ be defined as follows:
$$
\bar{g}(v,s_*)=\prod^d_{i=1}(2[v_i=*]+[v_i\neq*]).
$$

The proofs of the last two results in this section can be found in \cite{Kis}.
\begin{lemat}
\label{2dd}
Let $V\subset (*S)^d$. The set $V$ is a partition code if and only if $\sum_{v\in V} \bar{g}(v,s_*)= 2^d$.
\hfill{$\square$} 
\end{lemat}

\begin{wn}
\label{uV}
Let $V\subset S^d$ be a polybox code and let $u\in S^d$. For every $v\in V$ let $\bar{v}\in (*S)^d$ be defined  in the following way: if  $v_i\neq u_i$, then $\bar{v}_i=v_i$, and if  $v_i=u_i$, then $\bar{v}_i=*$. Let $\breve{u}\cap \breve{v}\neq\emptyset$ for every $v\in V$. If $u\sqsubseteq V$, then $\bar{V}=\{\bar{v}:v\in V\}$ is a partition code. 
\hfill{$\square$}
\end{wn}

\section{Small polybox codes}
To show that Keller's conjecture is true in dimension seven for a cube tiling $[0,1)^7+T$ for which $|L(T,x,i)|\geq 6$, it was sufficient to prove the following theorem (\cite{Kis}):

\begin{tw}
\label{12}
If $V,W\subset S^d, d\geq 4$, are disjoint and equivalent polybox codes without twin pairs, then $|V|\geq 12$.
\end{tw} 

To  show that the conjecture is true in dimension seven for a cube tiling $[0,1)^7+T$ with $|L(T,x,i)|=5$, we have to know the structure of all twin pairs free disjoint and equivalent polybox codes $V$ and $W$, with 12 words each, in dimensions four, five and six. To find this structure we need to know the structure of some polybox codes having a few words.

\subsection{Geometry of dichotomous boxes}
Before reading the proofs it is worth  paying attention to the basic aspects of the geometry of boxes that form two realizations $E(V)$ and $E(W)$ of equivalent polybox codes $V$ and $W$ (i.e., when $\bigcup E(V)=\bigcup E(W)$). 

\bigskip
{\center
\includegraphics[width=12cm]{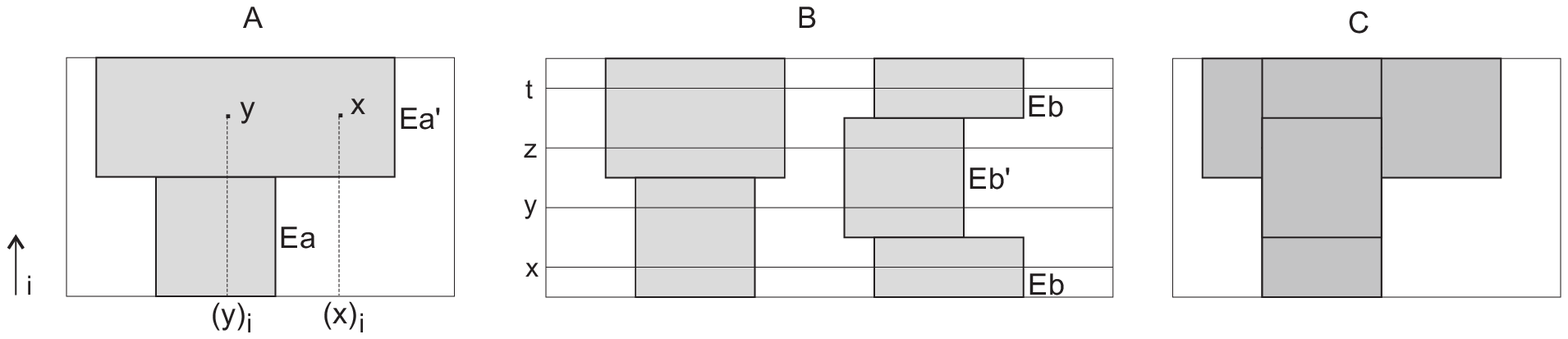}\\
}

\smallskip
\noindent{\footnotesize 
Fig. 3. Realizations $E(V)$ (A), $E(U)$ (B) and $E(W)$ (C), schematically, where $V=V^{i,a}\cup V^{i,a'}$,  $U=U^{i,a}\cup U^{i,a'}\cup U^{i,b}\cup U^{i,b'}$ and $W=W^{i,a'}\cup W^{i,b}\cup W^{i,b'}$. We assume that $\bigcup E(V)=\bigcup E(W)$.    
  }

\medskip
\noindent
(\textbf{P}): {\bf Projections}. Let $V^{i,a}$ and $V^{i,a'}$ be non-empty sets, and let $x\in \bigcup E(V^{i,a'})$ be such that $(x)_i\not\in \bigcup E((V^{i,a})_i$ (Figure 3A). Since words in $V$ and $W$ are dichotomous, $\bigcup E(V)=\bigcup E(W)$ and (\ref{dkostki}), we have $x\in \breve{w}$, where $w\in W$ is such that $w_i=a'$.

\smallskip
\noindent
(\textbf{S}): {\bf Slices}. By (\ref{dkostki}) for every $r\in El\cap Es$, $l\not\in \{s,s'\}$ the set $\pi^i_r=ES\times \cdots \times ES\times \{r\}\times ES\times \cdots \times ES$, where $\{r\}$ stands at the $i$-th position, slices the sets $\bigcup E(U^{i,l})$ and $\bigcup E(U^{i,s})$ simultaneously (Figure 3B, where $r\in \{x,y,z,t\}$).
  
\smallskip
\noindent
(\textbf{V}): {\bf Volumes}. Let $|V^{i,a}|=n$ and $|V^{i,a'}|=m$, and let $n< m$. Since all boxes $\breve{u}$, $u\in S^d$, are of the same size and $n<m$, by ({\bf P}),
$|W^{i,a'}|\geq m-n$.

\smallskip
\noindent
(\textbf{C}): {\bf Cylinders}. Suppose that $V^{i,b}\cup V^{i,b'}=\emptyset$ and $W^{i,b}\cup W^{i,b'}\neq\emptyset$. Then $\bigcup E((W^{i,b}))_i)=\bigcup E((W^{i,b'}))_i)$, and hence the set $\bigcup E((W^{i,b}\cup W^{i,b'}))$ in an $i$-cylinder in $(ES)^d$ (compare Figure 3A and 3C). By (\ref{2d}),  the codes $(W^{i,b})_i$ and $(W^{i,b'})_i$ are equivalent.


\smallskip
\noindent
(\textbf{Co}): {\bf Coverings}. Suppose that $(v)_i=(w)_i$, where $v\in V^{i,l}$ and $w\in W^{i,s}$, $l\not\in \{s,s'\}$. Then $(w)_i\sqsubseteq (V^{i,l'})_i$ and $(v)_i\sqsubseteq (W^{i,s'})_i$ (Figure 3A and 3C).

\bigskip

\subsection{Small partition codes}
All realizations of partition and polybox codes which are considered in this and the next subsection, may be replaced by one in which boxes $K=K_1\times \cdots \times K_d\subset [0,1]^d$ are such that $K_i\in \{[0,1],[0,1/2),[1/2,1]\}$. All considered suits are illustrated based on this realization. Keeping this in mind, our task in this subsection can be summarized as follows: for a given positive integer $n$ find all possible partition $\ka F$ of $[0,1]^d$  into pairwise dichotomous boxes $K$ such that $|\ka F|=n$. 

If $v\in (*S)^d$, and $\sigma$ is a permutation of the set $[d]$, then $v_\sigma=v_{\sigma(1)}\cdots v_{\sigma(d)}$. For every $i\in [d]$ let $h_i:*S\rightarrow *S$ be a bijection such that $h_i(*)=*$ and $(h_i(l'))'=h_i(l)$ for every $l\in S$. We say that polybox codes $V,W\subset (*S)^d$ are {\it isomorphic} if there are $\sigma$ and $h_1,\ldots ,h_d$ such that $W=\{h_1(v_{\sigma(1)})\cdots h_d(v_{\sigma(d)}): v\in V\}$. 

\begin{lemat}
\label{34}
Let $V\subset (*S)^d$ be a partition code. 

If $|V|=3$, then
\begin{equation}
\label{p3}
(V)_{A^c}=\{l_1*,l_1'l_2,l_1'l_2'\},
\end{equation}
where $l_1,l_2\in S$, $A=\{i_1<i_2\}\subseteq [d]$ and $(V)_A=\{*\cdots *\}\subset (*S)^{d-2}$.

If $|V|=4$ and $V$ contains only one twin pair, then
\begin{equation}
\label{p4}
(V)_{A^c}=\{l_1l_2l_3',l_1'l_2l'_3,*l'_2l'_3,**l_3\}
\end{equation}
where $l_1,l_2,l_3\in S$, $A=\{i_1<i_2<i_3\}\subseteq [d]$ and $(V)_A=\{*\cdots *\}\subset (*S)^{d-3}$.

If $|V|=5$ and $V$ does not contain a twin pair, then
\begin{equation}
\label{p5}
(V)_{A^c}=\{l_1l_2l_3,l_1'l_2'l_3',*l'_2l_3,l_1*l'_3,l'_1l_2*\},
\end{equation}
where  $l_1,l_2,l_3\in S$, $A=\{i_1<i_2<i_3\}\subseteq [d]$ and $(V)_A=\{*\cdots *\}\subset (*S)^{d-3}$.

If $|V|=6$ and $V$ does not contain a twin pair, then
\begin{equation}
\label{p6}
(V)_{A^c}=\{***l_4,\;l_1l_2l_3l'_4,\;l_1'l_2'l_3'l'_4,\;*l'_2l_3l'_4,\;l_1*l'_3l'_4,\;l'_1l_2*l'_4\},
\end{equation}
where  $l_1,l_2,l_3,l_4\in S$, $A=\{i_1<i_2<i_3<i_4\}\subseteq [d]$ and $(V)_A=\{*\cdots *\}\subset (*S)^{d-4}$.

The above partition codes are given up to an isomorphism.
\end{lemat}

{\center
\includegraphics[width=8cm]{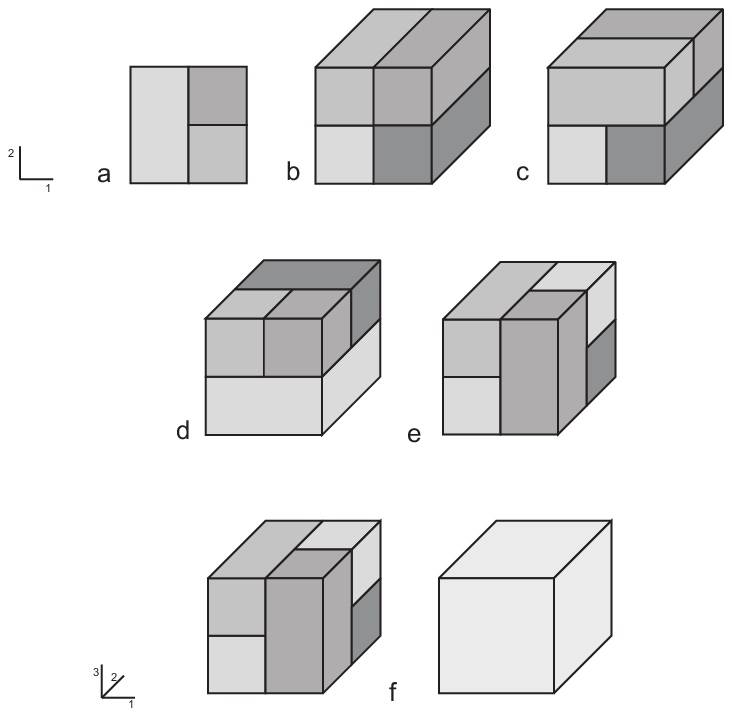}\\
}

\medskip
\noindent{\footnotesize Fig. 4. Figure a: the realization $f((V)_{A^c})$ of the partition code $(V)_{A^c}=\{l_1*,l_1'l_2,l_1'l_2'\}$ in the $2$-box $X=[0,1]^2$, where $f_i(l_i)=[0,1/2), f_i(*)=[0,1]$ for $i=1,2$. Figure b and c: examples of partitions of a $3$-box $X=[0,1]^3$ into four pairwise dichotomous boxes with more than one twin pair. Figure d:  the realization $f((V)_{A^c})$ in $X=[0,1]^3$, where $f_3(l_3)=[0,1/2)$ and $f_3(*)=[0,1]$,  of the partition code $(V)_{A^c}=\{l_1l_2l'_3,l'_1l_2l'_3,*l'_2l'_3,**l_3\}$ with one twin pair. Figure e: the realization $f((V)_{A^c})$ in $X=[0,1]^3$ of the partition code $(V)_{A^c}=\{l_1l_2l_3,l_1'l_2'l_3',*l'_2l_3,l_1*l'_3,l'_1l_2*\}$. Figure f: the realization $f((V)_{A^c})$ in the $4$-box $X=[0,1]^4$, where $f_4(l_4)=[0,1/2)$, $f_4(*)=[0,1]$ of the the code $(V)_{A^c}=\{***l'_4,\;l_1l_2l_3l_4,\;l_1'l_2'l_3'l_4,\;*l'_2l_3l_4,\;l_1*l'_3l_4,\;l'_1l_2*l_4\}$. We consider the two halves of the $4$-box $X=[0,1]^4$: $[0,1]^3\times [0,1/2)$ (on the left) and $[0,1]^3\times [1/2,1]$ (on the right), the fourth axis is omitted, in which we see the realizations of the codes  $(\{l_1l_2l_3l_4,\;l_1'l_2'l_3'l_4,\;*l'_2l_3l_4,\;l_1*l'_3l_4,\;l'_1l_2*l_4\})_{i_4}$ (on the left) and $(\{***l'_4\})_{i_4}$ (on the right).     
  }

\smallskip
\proof Let $V=\{v^1,v^2,v^3\}$. By Lemma \ref{2dd}, $\sum_{v\in V}\bar{g}(v,s_*)=2^d$, and thus $ \bar{g}(v^1,s_*)=2^{d-1}$, $\bar{g}(v^2,s_*)=\bar{g}(v^3,s_*)=2^{d-2}.$ Let $i_1\in [d]$ be such that $v^1_{i_1}\neq *$. Then $v^2_{i_1}=v^3_{i_1}=(v^1_{i_1})'$. The words $v^2,v^3$ are dichotomous, and therefore $v^2_{i_2}=(v^3_{i_2})'$, $v^2_{i_2}\neq *$, for some $i_2\in [d]\setminus \{i_1\}$ (Figure 4a). (Clearly, we can assume that $i_1<i_2.$) Obviously, $(V)_A=\{*\cdots *\}\subset (*S)^{d-2}$, where $A=\{i_1,i_2\}$.

Let now $V=\{v^1,v^2,v^3,v^4\}$. There are two solutions of the equation $\sum_{v\in V}\bar{g}(v,s_*)=2^d$: 
$\bar{g}(v^i,s_*)=2^{d-2}$ for every $i\in [4]$ and $ \bar{g}(v^1,s_*)=2^{d-1}$, $\bar{g}(v^2,s_*)=2^{d-2}$, $\bar{g}(v^3,s_*)=\bar{g}(v^4,s_*)=2^{d-3}.$ Since the words are pairwise dichotomous, it can be easily checked that in both cases there are $i_1,i_2,i_3\in [d]$, $i_1<i_2<i_3$, such that $v_i=*$ for every $v\in V$ and $i\in [d]\setminus\{i_1,i_2,i_3\}$. Thus, we have to determine all 
partitions of a $3$-dimensional box into four pairwise dichotomous boxes with only one twin pair. It is easy to see that the first solution corresponds to partitions with more than one twin pair (examples of such partitions are presented in Figure 4b and 4c). The second solution corresponds to  partition codes with one twin pair (Figure 4d).


The proofs of (\ref{p5}) and (\ref{p6}) (Figure 4e and  4f) can be found in \cite{Kis}.
\hfill{$\square$}

\subsection{Small polybox codes}
   
Now our goal is as follows: for given positive integers $n,m$ find all polyboxes $F\subset [0,1]^d$ having two suits $\ka F$ and $\ka G$ consisting of boxes $K$ and such that  $\ka F$ and $\ka G$ are disjoint, twin pairs free and $|\ka F|=n$, $|\ka G|=m$.

\begin{lemat}
\label{24}
Let $V,W\subset (*S)^d$  be disjoint and equivalent polybox codes without twin pairs such that $|V|\in \{2,3\}$ and $|W|\in \{2,3,4\}$, and let $l_1,l_2,l_3,l_4\in S$. 
\vspace{-20pt}
\begin{center}
\begin{displaymath}
\begin{tabular}{|l|r|}
\hline
\multicolumn{2}{|c|}{If $|V|=|W|=2$ and $d\geq 2$, then} \\ 
\hline
 $(V)_{A^c}=\{*l_2,\;l_1l_2'\}$ & $(W)_{A^c}=\{l'_1l_2,\;l_1*\}$ \\ \cline{1-2}
  \multicolumn{2}{|c|}{where $A=\{i_1,i_2\}$, $(V)_{A}=(W)_{A}=\{(p)_{A}\}$ for some $p\in (*S)^d$ }\\
\hline
\end{tabular}
\end{displaymath}
\end{center}
\vspace{-20pt}
\begin{center}
\begin{displaymath}
\begin{tabular}{|l|r|}
\hline
\multicolumn{2}{|c|}{If $|V|=2$, $|W|=3$ and $d\geq 3$, then} \\ 
\hline
 \multirow{1}{*}{$(V)_{A^c}=\{**l_3,\;l_1l_2l_3'\}$} & $(W)_{A^c}=\{l_1l_2*,\;l_1'l_2l_3,*l_2'l_3\}$ \\ \cline{1-2}
  $(V)_{A^c}=\{*l_2l_3,\;l_1*l_3'\}$ & $(W)_{A^c}=\{l_1l_2*,\;l_1'l_2l_3,\;l_1l_2'l_3'\}$ \\ \cline{1-2}
 \multicolumn{2}{|c|}{where $A=\{i_1,i_2,i_3\}$, $(V)_{A}=(W)_{A}=\{(p)_{A}\}$ for some $p\in (*S)^d$}\\
\hline
\end{tabular}
\end{displaymath}
\end{center}
\vspace{-20pt}
\begin{center}
\begin{displaymath}
\begin{tabular}{|l|r|}
\hline
\multicolumn{2}{|c|}{If  $|V|=2$, $|W|=4$ and $d\geq 3$, then} \\ 
\hline 
 \multirow{1}{*}{$(V)_{A^c}=\{***l_4,\;l_1l_2l_3l_4'\}$} & $(W)_{A^c}=\{l_1l_2l_3*,\;l_1'l_2l_3l_4,\;*l_2'l_3l_4,\;**l_3'l_4\}$ \\ \cline{1-2}
 \multirow{1}{*}{$(V)_{A^c}=\{**l_3l_4,\;l_1l_2*l'_4\}$} & $(W)_{A^c}=\{l_1l_2l_3*,\;l_1'l_2l_3l_4,\;*l_2'l_3l_4,\;l_1l_2l_3'l'_4\}$ \\ \cline{2-2}
 \cline{1-2}
 {$(V)_{B^c}=\{**l_3,\;l_1*l_3'\}$} & $(W)_{B^c}=\{l_1l_2*,\;l_1'l_2l_3,\;*l_2'l_3,\;l_1l_2'l_3'\}$ \\ \cline{1-2}
  \multicolumn{2}{|c|}{ where $A=\{i_1,i_2,i_3,i_4\}$, $(V)_{A}=(W)_{A}=\{(p)_{A}\}$ for some $p\in (*S)^d$,}  \\ 
  \multicolumn{2}{|c|}{ and $B=\{i_1,i_2,i_3\}$, $(V)_{B}=(W)_{B}=\{(p)_{B}\}$ for some $p\in (*S)^d$}\\ 
      \hline 
\end{tabular}
\end{displaymath}
\end{center}
\vspace{-20pt}
\begin{center}
\begin{displaymath}
\begin{tabular}{|l|r|}
\hline
\multicolumn{2}{|c|}{If $|V|=|W|=3$ and $d=3$, then} \\ 
\hline
 $V=\{*l_2l_3,l_1*l_3',l_1'l_2'*\}$ & $W=\{*l_2'l_3',l_1'*l_3,l_1l_2*\}$ \\ \cline{1-2}
 \multirow{1}{*}{$V=\{l_1**,l_1'l_2'*,l_1'l_2l_3\}$ } 
 &  $W=\{**l_3,*l_2'l_3',l_1l_2l_3'\}$ \\ \cline{1-2}
 $V=\{l_1l_2*,*l_2'l_3',l_1'l_2'l_3\}$ & $W=\{l_1'l_2'*,l_1*l_3',l_1l_2l_3\}$ \\ \cline{1-2}
  \hline
\end{tabular}
\end{displaymath}
\end{center}

The above polybox codes are given up to an isomorphism.\\
\end{lemat}

{\center
\includegraphics[width=12cm]{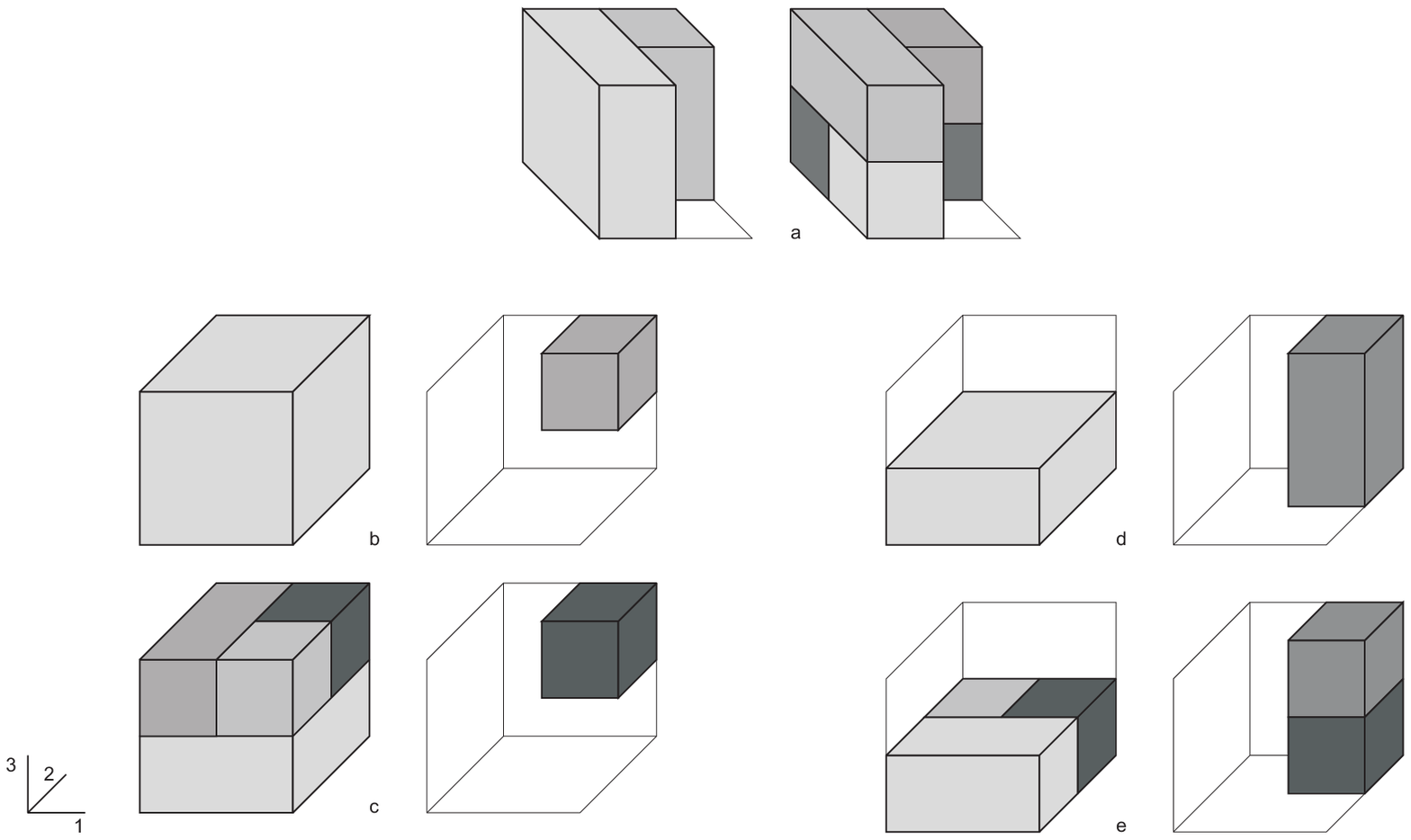}\\
}
\noindent{\footnotesize Fig. 5.  Figure a: the realizations $f((V)_{A^c})$ and $f((W)_{A^c})$ of the codes $(V)_{A^c}=\{l_1**,l_1'l_2'*\}$ (on the left) and $(W)_{A^c}=\{l_1l_2l_3,l_1'l_2'l_3',*l'_2l_3,l_1*l_3'\}$ (on the right) in the $3$-box $X=[0,1]^3$, where $f_i(l_i)=[0,1/2)$ and $f_i(*)=[0,1]$ for $i=1,2,3$. In a four dimensional case we consider the two halves of the $4$-box $X=[0,1]^4$: $[0,1]^3\times [0,1/2)$ (always on the left) and $[0,1]^3\times [1/2,1]$ (always on the right). Clearly, the fourth axis is omitted. Figure b and c: the realizations $(f((V)_{A^c}))_{i_4}$ and $(f((W)_{A^c}))_{i_4}$ of the codes $(V)_{A^c}=\{***l_4,l'_1l'_2l'_3l_4'\}$ and $(W)_{A^c}=\{l'_1l'_2l'_3*,l_1'l_2l_3'l_4,**l_3l_4,l_1*l_3'l_4\}$  in the $4$-box $X=[0,1]^4$, where $f_4(l_4)=[0,1/2)$ and $f_4(*)=[0,1]$. The two darkest boxes in Figure c are in fact one box which is a realization of the word $l'_1l'_2l'_3*$. Figure d and e: the realizations $(f((V)_{A^c}))_{i_4}$ and $(f((W)_{A^c}))_{i_4}$ of the codes $(V)_{A^c}=\{**l_3l_4,l'_1l'_2*l_4'\}$ and $(W)_{A^c}=\{l'_1l'_2l_3*,l_1l'_2l_3l_4,*l_2l_3l_4,l'_1l_2'l'_3l'_4\}$ in $X=[0,1]^4$. Similarly like above, the two darkest boxes in Figure e are one box which is a realization of the word $l'_1l'_2l_3*$.    
  }

\medskip
\noindent
{\it Proof of the case $|V|=2, |W|=4$.} Let $V=\{u,v\}$ and $W=\{w,p,q,r\}$.

We can assume that $\breve{u}\cap \breve{w}\neq\emptyset$ and $\breve{v}\cap \breve{w}\neq\emptyset$. Then there is $i\in [d]$ such that $(\breve{v})_i\cap (\breve{u})_i\neq\emptyset$ and $(\breve{w})_i\subseteq (\breve{v})_i\cap (\breve{u})_i$ (compare Figure 4A and C). We will show that $w_i=*$. Suppose this is not true. Since $(\breve{v})_i\cap (\breve{u})_i\neq\emptyset$ and $u,v$ are dichotomous, we have $v_i=u_i'$, $u_i\neq *$. Then $w_i\not\in \{u_i,u_i'\}$, and, by (\ref{dkostki}), we can choose $x\in \breve{u}\setminus \breve{w}$ and   $y\in \breve{v}\setminus \breve{w}$ such that $(x)_i=(y)_i$ and $(x)_i\in (\breve{w})_i$. The words in $W$ are pairwise dichotomous, and thus there is a word in $W$, say $p$, such that $x,y\in \breve{p}$. Note that $p_i=w_i'$ and consequently, $(\breve{p})_i\subseteq (\breve{v})_i\cap (\breve{u})_i$. 

Moreover, $(\breve{u})_i\setminus (\breve{v})_i\cup (\breve{v})_i\setminus (\breve{u})_i\neq\emptyset$, for otherwise $u$ and $v$ would be a twin pair. Let $(\breve{v})_i\setminus (\breve{u})_i\neq\emptyset$ and take $z\in \breve{v}$ such that $(z)_i\in (\breve{v})_i\setminus (\breve{u})_i$. Clearly, $z\not\in \breve{w}\cup \breve{p}$, and thus, $z\in \breve{q}$. Then, by (\textbf{P}) in Section 3.1, $q_i=v_i$. Since $w$ and $p$ are not a twin pair,  $(\breve{w})_i\setminus  (\breve{p})_i\cup (\breve{p})_i\setminus  (\breve{w})_i\neq\emptyset$. Assume without loss of generality that $(\breve{p})_i\setminus  (\breve{w})_i\neq\emptyset$ and choose $z^1\in \breve{u}\setminus \breve{p}$ such that $(z^1)_i\in (\breve{p})_i\setminus  (\breve{w})_i$. Then $z^1\in \breve{r}$, and since $p$ and $r$ are dichotomous, we have $r_i=p_i'=w_i$. Now it can be easily seen that $(\breve{w})_i\cup (\breve{r})_i=(\breve{u})_i$, which implies that $w$ and $r$ form a twin pair, a contradiction.
This completes the proof that $w_i=*$. 

We now show that exactly one box from the set $E(W)$  has nonempty intersection with both boxes $\breve{u}$ and $\breve{v}$. Assume on the contrary that there are exactly two boxes in $E(W)$, say $\breve{w}$ and $\breve{p}$, having nonempty intersections with $\breve{u}$ and $\breve{v}$ simultaneously. Then, as we have just shown, $w_i=p_i=*$ and $q_i,r_i\in \{u_i,u_i'\}$. If $(\breve{u})_i=(\breve{w})_i\cup (\breve{p})_i$ or  $(\breve{v})_i=(\breve{w})_i\cup (\breve{p})_i$, then $(\breve{w})_i$ and $(\breve{p})_i$ are a twin pair, and consequently, $w$ and $p$ are a twin pair, which is a contradiction. Therefore, $(\breve{u})_i=(\breve{w})_i\cup (\breve{p})_i\cup (\breve{q})_i$ and  $(\breve{v})_i=(\breve{w})_i\cup (\breve{p})_i\cup (\breve{r})_i$ and $q_i=r_i'$. By (\ref{p3}), there are  twin pairs in the sets of boxes $\{(\breve{w})_i,(\breve{p})_i,(\breve{q})_i\}$ and  $\{(\breve{w})_i,(\breve{p})_i,(\breve{r})_i\}$.
Since $w_i=p_i$, the boxes $(\breve{w})_i$ and $(\breve{p})_i$ cannot form a twin pair which means that the set $(\breve{w})_i\cup (\breve{p})_i$ is not a box. But $(\breve{w})_i\cup (\breve{p})_i\cup (\breve{q})_i$ and  $(\breve{w})_i\cup (\breve{p})_i\cup (\breve{r})_i$ are boxes. Therefore, $(\breve{q})_i=(\breve{r})_i$ and consequently, $q$ and $r$ are a twin pair, a contradiction.

If $w_i=p_i=q_i=*$ and $r_i\neq *$, then  $(\breve{w})_i\cup (\breve{p})_i\cup (\breve{q})_i=(u)_i$ or  $(\breve{w})_i\cup (\breve{p})_i\cup (\breve{q})_i=(v)_i$. By (\ref{p3}), two of the three boxes $(\breve{w})_i,(\breve{p})_i,(\breve{q})_i$ are a twin pair, and therefore there is a twin pair among the words $w,p,q$, which is impossible. 

Similarly, if  $w_i=p_i=q_i=r_i=*$, then $(\breve{w})_i\cup (\breve{p})_i\cup (\breve{q})_i\cup (\breve{r})_i=(\breve{u})_i$ and therefore, by (\ref{p4}) and the proof of Lemma \ref{34} (the case $|V|=4$), there is a twin pair in the set $\{w,p,q,r\}$, which is a contradiction.  
 
We have shown that $w_i=*$, $u_i=v_i'$ and $p_i,q_i,r_i\in \{u_i,u_i'\}$.

Since $(\breve{w})_i\subseteq (\breve{u})_i$, for every $j\in [d]\setminus \{i\}$ we have $Ew_j\subseteq Eu_j$, and by (\ref{dkostki}), if $Ew_j\neq Eu_j$, then $w_j\in S$ and $u_j=*$. 

Now two cases may occur: 
$$
\breve{u}=\breve{u}\cap \breve{w}\cup \breve{p}\cup \breve{q}\cup \breve{r}\;\; {\rm and}\;\; \breve{v}=\breve{v}\cap \breve{w} \;\;\;\;\;({\rm Figure\; 5b,c})
$$
or 
$$
\breve{u}=\breve{u}\cap \breve{w}\cup \breve{p}\cup \breve{q}\;\; {\rm and}\;\; \breve{v}=\breve{v}\cap \breve{w}\cup \breve{r} \;\;\;\;\;({\rm Figure\; 5a,d,e}).
$$ 

In the first case the $d$-box $\breve{u}$ is partitioned into four pairwise dichotomous boxes, and thus the structure of this partition is given by (\ref{p4}) which contains exactly one twin pair. Hence, the partition $\{\breve{u}\cap \breve{w}, \breve{p}, \breve{q}, \breve{r}\}$ contains one twin pair, and $W$ does not contain a twin pair. Therefore, the box $\breve{u}\cap \breve{w}$ must be  one of the twins. We may assume that $\breve{p}$ is the second one. Thus, there are $i_1,i_2,i_3\in [d]$, $i_1<i_2<i_3$ and the letters  $l_1,l_2,l_3\in S$ such that (we assume without loss of generality that $i_3<i$ 
): 
$(u)_{A^c}=***l_4, (w)_{A^c}=l_1l_2l_3*$, and therefore
$(\breve{u}\cap \breve{w})_{A^c}=El_1\times El_2\times El_3\times El_4$, where $l_4=u_i$, $A=\{i_1<i_2<i_3<i\}$ and  $(p)_{A^c}$ has one of the forms: $l_1'l_2l_3l_4$, $l_1l'_2l_3l_4$ or $l_1l_2l'_3l_4$. We consider the first case as in the rest of them we obtain isomorphic forms. Let  $(p)_{A^c}=l_1'l_2l_3l_4$. By (\ref{p4}), $(q)_{A^c}=*l_2'l_3l_4$ and $(r)_{A^c}=**l_3'l_4$ or $(q)_{A^c}=*l_2l'_3l_4$ and   $(r)_{A^c}=*l_2'*l_4$. Since $(w)_i=(v)_i$, we have $(v)_{A^c}=l_1l_2l_2l_4'$. By (\ref{p4}), $(u)_A=(p)_A=(q)_A=(r)_A$ and, since $(\breve{w})_i\subset (\breve{u})_i\cap (\breve{v})_i$, $(u)_A=(w)_A=(v)_A$.   

In the second case there are two possibilities: $(\breve{u})_i\cap (\breve{r})_i\neq\emptyset$ or $(\breve{u})_i\cap (\breve{r})_i=\emptyset$. 
Since now the $d$-box $\breve{u}$ is divided into three pairwise dichotomous boxes, the structure of the partition $\{\breve{u}\cap \breve{w}, \breve{p}, \breve{q}\}$ is given by (\ref{p3}). Clearly, as above, we may assume that the boxes  $\breve{u}\cap \breve{w}$ and $\breve{p}$ are the only twin pair in this partition. Thus, there are $i_1,i_2\in [d]\setminus\{i\}$, $i_1<i_2$, and letters $l_1,l_2\in S$ such that  such that $(u)_{B^c}=**l_3$, $(w)_{B^c}=l_1l_2*$ and $(u)_B=(w)_B$, where $l_3=u_i$ and $B=\{i_1<i_2<i_3\}$, $i=i_3$. Furthermore, $(p)_{B^c}=l_1'l_2l_3$ and $(q)_{B^c}=*l_2'l_3$ or  $(p)_{B^c}=l_1l'_2l_3$ and $(q)_{B^c}=l_1'*l_3$ (in this second case we obtain an isomorphic form). In both cases, $(p)_B=(q)_B=(u)_B$.

Let $(\breve{u})_i\cap (\breve{r})_i\neq\emptyset$. Since $(\breve{w}\cap \breve{v})_i$ and $(\breve{r})_i$ are a twin pair and $(u)_B=(w)_B$, there is  
$k\in \{i_1,i_2\}$ such that $w_k=r_k'$ and $(w)_{\{i,k\}}=(r)_{\{i,k\}}$. Then $(w)_B=(r)_B$. Taking  $k=i_1$, we obtain $(r)_{B^c}=l_1'l_2l_3'$. Thus, we have to exclude the case $(p)_{B^c}=l_1'l_2l_3$ and $(q)_{B^c}=*l_2'l_3$, for otherwise $p$ and $r$ are a twin pair (if we take $k=i_2$, then $(r)_{B^c}=l_1l_2'l_3'$ and the case $(p)_{B^c}=l_1l_2'l_3$ and $(q)_{B^c}=*l_2l_3$ has to be excluded). Since $\breve{v}=\breve{v}\cap \breve{w}\cup \breve{r}$, we have $(v)_{B^c}=*l_2l_3'$ and $(v)_{\{i_1,i\}}=(w)_{\{i_1,i\}}$, and then $(v)_B=(w)_B$ (if $k=i_2$, then $(v)_{B^c}=l_1*l_3'$). 

Summing up, in the case  $(\breve{u})_i\cap (\breve{r})_i\neq\emptyset$ we have: $(u)_{B^c}=**l_3$, $(v)_{B^c}=*l_2l_3'$ and  $(w)_{B^c}=l_1l_2*$, $(p)_{B^c}=l_1l'_2l_3$, $(q)_{B^c}=l_l'*l_3$ and $(r)_{B^c}=l_1'l_2l_3'$. Moreover, $(p)_B=(q)_B=(u)_B=(v)_B=(w)_B=(r)_B$.

Let now $(\breve{u})_i\cap (\breve{r})_i=\emptyset$. Since the boxes $(\breve{w}\cap \breve{v})_i$ and $(\breve{r})_i$ are a twin pair, if $w_{i_1}=r_{i_1}'$ or  $w_{i_2}=r_{i_2}'$, then $(\breve{u})_i\cap (\breve{r})_i\neq\emptyset$. Therefore there is exactly one $j\in B^c$ such that $w_j=r_j'$ and $w_j\neq *$. Assume without loss of generality that $j=i_4>i_3=i$. Then $(r)_{A^c}=l_1l_2l_3'l_4'$, where $l_4=u_{i_4}$($\neq *$ because $w_j=u_j$) and $A=\{i_1,i_2,i_3,i_4\}$. Thus, $(u)_{A^c}=**l_3l_4$, $(v)_{A^c}=l_1l_2l_3'*$ and $(w)_{A^c}=l_1l_2*l_4$. Clearly,  $(p)_{A^c}=l_1'l_2l_3l_4$ and $(q)_{A^c}=*l_2'l_3l_4$ (if   $(p)_{A^c}=l_1l'_2l_3l_4$ and $(q)_{A^c}=l_l'*l_3l_4$ we get an isomorphic form of $W$). Since $\breve{u}=\breve{p}\cup \breve{q}\cup \breve{w}\cap \breve{u}$, we have $(p)_A=(q)_A=(u)_A=(w)_A$. Similarly, since $\breve{v}=\breve{r}\cup \breve{w}\cap \breve{v}$, we have $(w)_A=(r)_A=(v)_A$.

Permuting the letters at the third and the fourth position in every word in $V$ and $W$ we get the form as it is in the lemma.
\hfill{$\square$}

\medskip
{\center
\includegraphics[width=11cm]{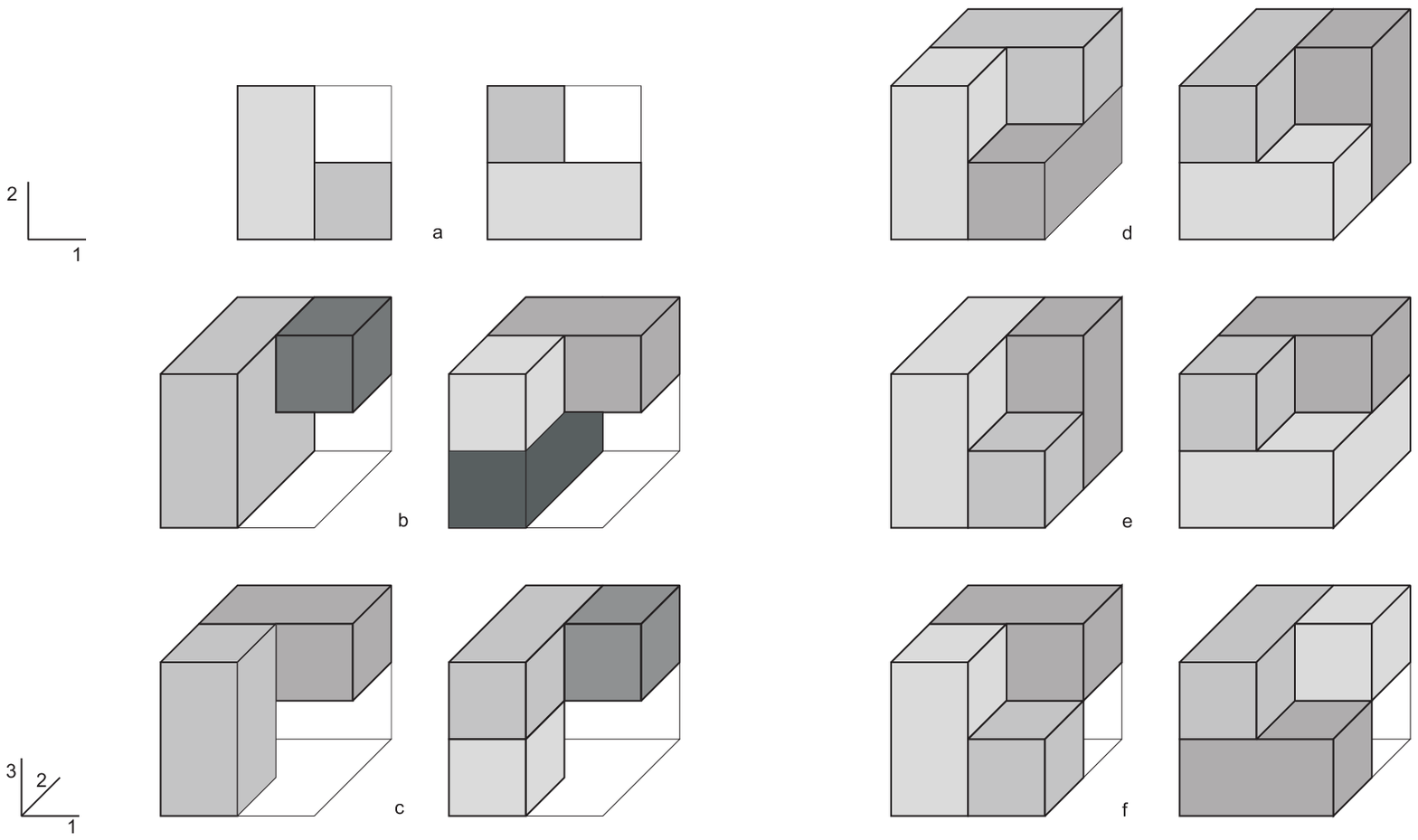}\\
}

\bigskip
\noindent{\footnotesize Fig. 6.  Figure a: the realizations $f((V)_{A^c})$ and $f((W)_{A^c})$ of the codes $(V)_{A^c}=\{l_1*,l_1'l_2\}$ (on the left) and $(W)_{A^c}=\{*l_2,l_1l_2'\}$ (on the right) in the $2$-box $X=[0,1]^2$, where $f_i(l_i)=[0,1/2)$ and $f_i(*)=[0,1]$ for $i=1,2$. Figure b: the realizations $f((V)_{A^c})$ and $f((W)_{A^c})$ of the codes $(V)_{A^c}=\{l_1**,l_1'l_2'l_3'\}$ and $(W)_{A^c}=\{l_1*l_3,l_1l_2l_3',*l_2'l_3'\}$ in the $3$-box $X=[0,1]^3$, where $f_3(l_3)=[0,1/2)$ and $f_3(*)=[0,1]$.  Figure c: the realizations $f((V)_{A^c})$ and $f((W)_{A^c})$ of the codes $(V)_{A^c}=\{l_1l_2*,*l_2'l_3'\}$ and $(W)_{A^c}=\{l_1l_2l_3,l_1*l_3',l_1'l_2'l_3'\}$. Figure d: the realizations $f(V)$ and $f(W)$ of the codes $V=\{l_1l_2*,l_1'*l_3,*l_2'l_3'\}$ (on the left) and $W=\{*l_2l_3, l_1*l_3', l_1'l_2'*\}$ (on the right).  Figure e: the realizations $f(V)$ and $f(W)$ of the codes $V=\{l_1**,l_1'l_2l_3,l_1'l_2'*\}$ and $W=\{**l_3,l_1l_2l_3',*l_2'l_3'\}$.  Figure f: the realizations $f(V)$ and $f(W)$ of the codes $V=\{l_1l_2*,l_1'l_2l_3,*l_2'l_3'\}$ and $W=\{*l_2l_3,l_1*l_3',l_1'l_2'l_3'\}$.    
  }

\medskip  
\noindent
{\it The sketches of the proofs of the rest of the cases of Lemma \ref{24}.} 

Let $|V|=2,|W|=3$, and let $V=\{v,u\}$ and $W=\{w,p,q\}$. In the same way as above we show that there is exactly one word in $W$, say $w$, and there is $i\in [d]$ such that $w_i=*$, $u_i=v_i'$, $(w)_i\subseteq (u)_i\cap (v)_i$ and $p_i,q_i\in \{u_i,u_i'\}$. 

If $\breve{u}=\breve{u}\cap \breve{w}\cup \breve{p}\cup \breve{q}$ and $\breve{v}=\breve{v}\cap \breve{w}$, then the structure the partition $\{\breve{u}\cap \breve{w},\breve{p},\breve{q}\}$ of the $d$-box $\breve{u}$ is given by (\ref{p3}), and $(\breve{v})_i=(\breve{w})_i$. This case is illustrated in Figure 6b.

If  $\breve{u}=\breve{u}\cap \breve{w}\cup \breve{p}$ and $\breve{v}=\breve{v}\cap \breve{w}\cup \breve{q}$, then the boxes $\breve{u}\cap \breve{w}$, $\breve{p}$ are  a twin pair and $\breve{v}\cap \breve{w}$, $\breve{q}$ are a twin pair. Note that $(\breve{u})_i\cap (\breve{q})_i=\emptyset$, for otherwise $p$ and $q$ are a twin pair, which is impossible.  This case is illustrated in Figure 6c.

The proof of the case  $|V|=2,|W|=2$ (Figure 6a) can be found in \cite{Kis}.

Let  $|V|=\{v^1,v^2,v^3\}$ and $|W|=\{w^1,w^2,w^3\}$. 
Since in this case realizations are $3$-dimensional boxes, we will establish only the values $\bar{g}(v^i,s_*)$ and $\bar{g}(w^i,s_*)$ for $i=1,2,3$. Recall that $\bar{g}(v^i,s_*)=2^k$ if and only if the words $v^i$ contains $k$ stars and, by Lemma \ref{2dd}, $\sum_{i=1}^3\bar{g}(v^i,s_*)\leq 8$.

If $\bar{g}(v^1,s_*)=4$ and $\bar{g}(v^2,s_*)=2$, then $\bar{g}(v^3,s_*)=1$ because if   $\bar{g}(v^3,s_*)=2$, then $v^2$ and $v^3$ form a twin pair. This is easy to see that  $\bar{g}(w^1,s_*)=4$, $\bar{g}(w^2,s_*)=2$, $\bar{g}(w^3,s_*)=1$. This case is illustrated in Figure 6e. 

It can be easily verified that the case $\bar{g}(v^1,s_*)=4$, $\bar{g}(v^2,s_*)=1$ and $\bar{g}(v^3,s_*)=1$ is impossible. 

Let now  $\bar{g}(v^i,s_*)=2$ for $i=1,2,3$. This is obvious that $\bar{g}(w^i,s_*)=2$ for $i=1,2,3$ (Figure 6d).

Similarly, this is not hard to find the forms of $V$ and $W$ in the case  when $\bar{g}(v^1,s_*)=\bar{g}(v^2,s_*)=2$  and $\bar{g}(v^3,s_*)=1$ (Figure 6f). 

Finally, this is easy to check that the cases $\bar{g}(v^1,s_*)=2$, $\bar{g}(v^2,s_*)=\bar{g}(v^3,s_*)=1$ and $\bar{g}(v^i,s_*)=1$ for $i=1,2,3$ are impossible.
\hfill{$\square$}

\medskip 
To examine the structure of a polybox code $V$ we will consider the partitions (into pairwise dichotomous boxes) of a $d$-box $\breve{w}$ of the form  $\{\breve{v}\cap \breve{w}: v\in V\}$, where $w\in S^d$. The structures of these partitions are described by the above presented codes, where instead of the letter $*$ we put the adequate letter $w_i$. This process is shown in the following example.

\medskip
{\center
\includegraphics[width=12cm]{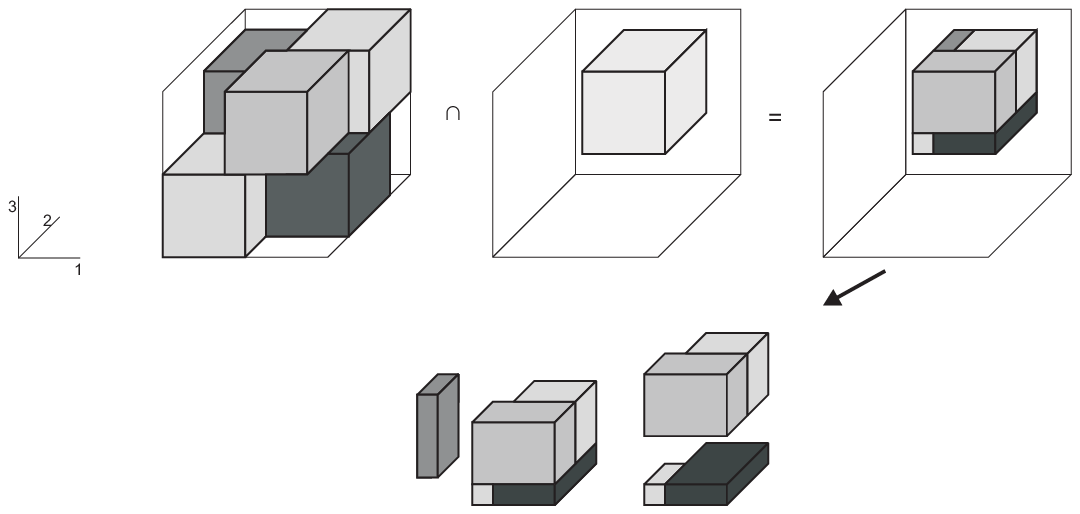}\\
}

\bigskip
\noindent{\footnotesize Fig. 7. On the top: the light box (in the middle) is contained in the sum of five pairwise dichotomous boxes (the boxes on the left). These boxes determine a partition of the light box into pairwise dichotomous boxes (the partition on the right). On the bottom: the boxes in this partition are arranged into $3$-cylinders.      
  }

\medskip
\begin{ex}
\label{p}
{\rm In Figure 7 the five boxes on the left are a realization of the polybox code $V=\{aaa,a'a'a',baa',a'ba, aa'b\}$, and the box in the middle is a realization of the word $w=bbb$. Since $w\sqsubseteq V$, we have $\breve{w}\subset \bigcup E(V)$. Thus, the $3$-box $\breve{w}$ is divided into pairwise dichotomous boxes $\breve{w}\cap \breve{v}$ for $v\in V$, and the set  $\bigcup (\{\breve{w}\cap \breve{v}: v\in Q\} \cup \{\breve{w}\cap \breve{v}: v\in P\})$, where  $P=\{v\in V^{3,a}: \breve{w}\cap \breve{v}\neq\emptyset\}=\{aaa,a'ba\}$ and $Q=\{v\in V^{3,a'}: \breve{w}\cap \breve{v}\neq\emptyset\}=\{a'a'a',baa'\}$, is an $3$-cylinder in the box $\breve{w}$. Therefore, $\bigcup \{(\breve{w}\cap \breve{v})_3: v\in Q\}=\bigcup \{(\breve{w}\cap \breve{v})_3: v\in P\}$. Now, the polybox $\bigcup \{(\breve{w}\cap \breve{v})_3: v\in Q\}$ is divided twice into pairwise dichotomous boxes without twin pairs. Since $|Q|=|P|=2$, we apply Lemma \ref{24} for the case $|V|=|W|=2$ to get the structure of $(Q)_3$ and $(P)_3$. Recall that in that case  we have $(V)_{A^c}=\{*l_2,\;l_1l_2'\}$ and $(W)_{A^c}=\{l'_1l_2,\;l_1*\}$, where $A=\{i_1,i_2\}$ and $(V)_{A}=(W)_{A}=\{(p)_{A}\}$. In our case we have $A=\{1,2\}$. Making in $(V)_{A^c}$ and $(W)_{A^c}$ the substitutions $l_1=a',l_2=a$ and $*=b$ we obtain $(Q)_3=\{ba,a'a'\}$ and $(P)_3=\{aa,a'b\}$.}
\end{ex}

\smallskip  
We now collect the above results in the forms in which they will be used later in the paper.

\begin{st}
\label{stt}
Let $V\subset S^d$ be a polybox code and let  $w\sqsubseteq V$, where $w\in S^d$ and $w\not\in V$.

\smallskip
(a) If the polybox code $V$ does not contain a twin pair, then it contains at least five words. The code $V$ has exactly five words and does not contain a twin pair if and only if it is of the form (\ref{p5}) of Lemma \ref{34}, where instead of $*$ at a position $j\in [d]$ we take $w_{j}$,  and $l_k\not\in \{w_{i_k},w_{i_k}'\}$ for $k\in \{1,2,3\}$. Moreover, $V$ is rigid. If $\breve{w}\cap \breve{v}\neq\emptyset$ for every $v\in V$ and the code $V$ has exactly six words and does not contain a twin pair, then it is of the form (\ref{p6}) of Lemma \ref{34}, where instead of $*$ at a position $j\in [d]$ we take $w_{j}$,  and $l_k\not\in \{w_{i_k},w_{i_k}'\}$ for $k\in \{1,2,3,4\}$.

(b) If $u\sqsubseteq V$, where $u\in S^d$ and $u\not\in V$, the code $V$ does not contain twin pairs and the words $w,u$ are dichotomous but they do not form a twin pair, then $|V|\geq 7$.

(c) Let  $\{\breve{w}\cap \breve{v}\neq\emptyset: v\in V^{i,l}\} \neq\emptyset$, where $w_i\not\in \{l,l'\}$. Then the set $\breve{w}\cap \bigcup E(V^{i,l})\cup \breve{w}\cap \bigcup E(V^{i,l'})$ is an $i$-cylinder in the $d$-box $\breve{w}$. Consequently, 
$$
\bigcup \{(\breve{w}\cap \breve{v})_i: v\in Q\}=\bigcup \{(\breve{w}\cap \breve{v})_i: v\in P\},
$$ where $Q=\{v\in V^{i,l}: \breve{w}\cap \breve{v}\neq\emptyset\}$ and $P=\{v\in V^{i,l'}: \breve{w}\cap \breve{v}\neq\emptyset\}$.

If $|P|=1$, $|Q|=5$  and $V$ does not contain a twin pair, then there is a set $A=\{i_1<i_2<i_3\}\subset [d]$ such that 
$$
(P)_{A^c}=\{w_{i_1}w_{i_2}w_{i_3}\}\;\; {\rm and}\;\;(Q)_{A^c}=\{l_1l_2l_3,l_1'l_2'l_3',w_{i_1}l'_2l_3,l_1w_{i_2}l'_3,l'_1l_2w_{i_3}\}, 
$$
where $l_k\not\in \{w_{i_k},w_{i_k}'\}$ for $k=1,2,3$ and $(P)_{A\cup \{i\}}=(Q)_{A\cup \{i\}}$

If  $(|P|,|Q|)\in \{(2,2),(2,3),(2,4),(3,3)\}$  and $V$ does not contain a twin pair, then the structure of $(P)_i$ and $(Q)_i$ is such as in Lemma \ref{24}, but in all those polybox codes we put $w_j$ instead of $*$, if the star appears at the $j$-th position, and $l_k\not\in \{w_{i_k},w_{i_k}'\}$ for $k\in \{1,2,3,4\}$.

(d) Let $P$ and $Q$ be such as in (c). If $|P|=1$ and $1\leq |Q|\leq 4$, then there is a twin pair in $V$. 


\end{st}

\medskip
\noindent
{\it Proof of (a)}. For $|V|=5$ it can be found in \cite{Kis}, and the case $|V|=6$ is proven in the same manner (see also \cite{Kis}).
\hfill{$\square$}

\medskip
\noindent 
{\it Proof of (b).} Let $W=\{v\in V: \breve{v}\cap \breve{w}\neq\emptyset\}$ and  $U=\{v\in V: \breve{v}\cap \breve{u}\neq\emptyset\}$. By (a), $|W|\geq 5$ and $|U|\geq 5$. 

Suppose that  $|W|=5$ and $|U|=5$. Again by (a), there is a set $A=\{i_1<i_2<i_3\}\subset [d]$ and letters $l_1,l_2,l_3\in S$, $l_j\not\in \{w_{i_j},w_{i_j}'\}$ for $i=1,2,3$, such that  
$$
(W)_{A^c}=\{l_1l_2l_3,l_1'l_2'l_3',w_{i_1}l'_2l_3, l_1w_{i_2}l'_3, l'_1l_2w_{i_3}\}
$$
and $(w)_A=(v)_A$ for every $v\in W$.

Clearly, if $|W\cap U|\leq 3$, then $|V|\geq 7$. 

Let $|W\cap U|=4$. Since $w$ and $u$ are dichotomous, there is $i\in [d]$ such that $w_i=u_i'$. If $i\in A^c$, then $W\cap U=\emptyset$ because $(w)_A=(v)_A$ for every $v\in W$. Therefore, $i\in A$. Observe that, for the rest two  $j\in A\setminus \{i\}$ we have $u_j\neq w_j'$ and $u_j\not\in \{l_j, l_j'\}$, for otherwise  $|W\cap U|<4$. Assume without loss of generality that $i=i_1$. Since the structure of $U$ is such  as predicted in (a), $|W\cap U|=4$ and $w_i=u_i'$, 
 it follows that
$$
(U)_{A^c}=\{l_1l_2l_3,l_1'l_2'l_3',w'_{i_1}l'_2l_3, l_1w_{i_2}l'_3, l_1'l_2w_{i_3}\},
$$
and $(w)_A=(u)_A$. Thus, $w'_{i_1}=u_{i_1}$, $w_{i_2}=u_{i_2}$ and $w_{i_3}=u_{i_3}$, and hence $w$ and $u$ are a twin pair, a contradiction.

If $W=U$, then $u=w$, which contradicts the assumption.

Suppose now that $|W|=6$ and assume on the contrary that $|V|=6$. Then $V=W$. By (\ref{p6}) in Lemma \ref{34} and (a), 
$$    
(W)_{A^c}=\{w_{i_1}w_{i_2}w_{i_3}l_4,\;l_1l_2l_3l'_4,\;l_1'l_2'l_3'l'_4,\;w_{i_1}l'_2l_3l'_4,\;l_1w_{i_2}l'_3l'_4,\;l'_1l_2w_{i_3}l'_4\},
$$
and $(w)_A=(v)_A$ for $v\in W$, where $A=\{i_1<i_2<i_3<i_4\}\subseteq [d]$ and $l_j\not\in \{w_{i_j}.w_{i_j}'\}$ for $j=1,2,3,4$. 

Let $v^1\in W$ be such that $(v^1)_{A^c}=w_{i_1}w_{i_2}w_{i_3}l_4$. We have $(w)_{i_4}=(v^1)_{i_4}$. Note that the structure of $(W\setminus\{v^1\})_{i_4}$ is such as in (a). In particular, $(w)_{i_4}$ is one and only word which is covered by $(W\setminus\{v^1\})_{i_4}$. Hence, $(u)_{i_4}=(w)_{i_4}$. Consequently, $w$ and $u$ are a twin pair, a contradiction.  
\hfill{$\square$}

\medskip
\noindent
{\it Proof of (c).} The set  $\breve{w}\cap \bigcup E(V^{i,l})\cup \breve{w}\cap \bigcup E(V^{i,l'})$ is an $i$-cylinder in the $d$-box $\breve{w}$ because $\{\breve{w}\cap \breve{v}:v\in V\}$ is a suit for $\breve{w}$ (Figure 7).  

Since $V$ does not contain twin pairs, the set of boxes $\{\breve{w}\cap \breve{v}: v\in V\}$ is a partition of the $d$-box $\breve{w}$ into pairwise dichotomous boxes which, by Lemma \ref{=c}, does not contain twin pairs (Figure 7). 

Since the set $\bigcup \{\breve{w}\cap \breve{v}: v\in  V^{i,l}\}\cup \bigcup \{\breve{w}\cap \breve{v}: v\in  V^{i,l'}\}$ is an $i$-cylinder in the box $\breve{w}$, we have $\bigcup  \{(\breve{w}\cap \breve{v})_i: v\in  V^{i,l}\}=\bigcup \{(\breve{w}\cap \breve{v})_i: v\in  V^{i,l'}\}$. 

We prove only the case  $|P|=1$, $|Q|=5$. The rest of the cases is proven in the very similar way (compare Example \ref{p}).

Let $P=\{u\}$. The $(d-1)$-box $(\breve{w}\cap \breve{u})_i$ is divided into five pairwise dichotomous boxes $\{(\breve{w}\cap \breve{v})_i: v\in Q\}$. Thus, $(\breve{w}\cap \breve{v})_i\subseteq (\breve{w}\cap \breve{u})_i$ for every $v\in Q$, and then $Ew_j\cap Ev_j\subseteq Ew_j\cap Eu_j$ for every $j\in [d]\setminus \{i\}$. It follows that, by (\ref{dkostki}), if $w_j\neq u_j$, then $v_j=u_j$.
Moreover, by Lemma \ref{=c}, the boxes of the partition $\{\breve{w}\cap \breve{v}:v\in Q\}$ do not form twin pairs. 
Therefore, a code of the partition   $\{(\breve{w}\cap \breve{v})_i: v\in Q\}$ of the box $(\breve{w}\cap \breve{u})_i$ is given by (\ref{p5}). Since for every $j\in A=\{i_1<i_2<i_3\}$ ($A$ is such as in (a)) there is $l\in S\setminus \{w_j\}$ such that $v_j=l$ and $q_j=l'$ for some $v,q\in Q$ and $Ew_j\cap Eu_j=Ew_j\cap Ev_j\cup Ew_j\cap Eq_j$, it must be, by (\ref{dkostki}), $w_j=u_j$ for every $j\in A$. Thus, $(P)_{A^c}=\{w_{i_1}w_{i_2}w_{i_3}\}$ and $(Q)_{A^c}=\{l_1l_2l_3,l_1'l_2'l_3',w_{i_1}l'_2l_3,l_1w_{i_2}l'_3,l'_1l_2w_{i_3}\}$.  Since, by (\ref{p5}), $Ew_j\cap Eu_j=Ew_j\cap Ev_j$ for every $v\in Q$ and $j\in [d]\setminus (A\cup \{i\})$, we have, by (\ref{dkostki}), $u_j=v_j$, and thus $(P)_{A\cup \{i\}}=(Q)_{A\cup \{i\}}$. In particular, $(u)_i\sqsubseteq (Q)_i$. 
\hfill{$\square$}

\medskip
\noindent
{\it Proof of (d).} Let $P=\{u\}$. The $(d-1)$-box $(\breve{w}\cap \breve{u})_i$ is divided into pairwise dichotomous boxes $\{(\breve{u}\cap \breve{v})_i:v\in  Q\}$. By (\ref{p3}), (\ref{p4}) and the proof of Lemma \ref{34}, this partition contains a twin pair. Using similar arguments as in (c) we show that there is a twin pair in $Q$.
\hfill{$\square$}

\section{The structure of equivalent polyboxes codes with 12 words: necessary conditions }
In this section we determine  necessary conditions which have to be fulfilled by disjoint and equivalent twin pairs free polyboxes codes $V$ and $W$ having 12 words each. This conditions will serve us to establish the initial configurations of words for the computations.   

In \cite{Kis} we defined a graph on a polybox code $V$. Recall that, a pair of words $v,u\in S^d$ such that $v_i\not\in \{u_i,u_i'\}$ for some $i\in [d]$ and  $(u)_{i}$ and $(v)_{i}$ are a twin pair is called an $i$-{\it siblings} (in Figure 7 the two upper boxes on the left upper picture are $1$-siblings). Let $V\subset S^d$ be a polybox code.  {\it A graph of siblings in V} is a graph $G=(V,\ka E)$ 
in which two vertices $u,v\in V$ are adjacent if they are $i$-siblings for some $i\in [d]$. We colour each edge in $\ka E$ with the colours from the set $[d]$: an edge $e\in \ka E$ has a colour $i\in [d]$ if its ends are $i$-siblings. The graph $G$ is simple and, if $V$ does not contain a twin pair, $d(v)\leq d$ for every $v\in V$, where $d(v)$ denotes the number of neighbors of $v$. Observe that the graph $G$ does not contain triangles.

In \cite{Kis} we proved the following two lemmas.

\begin{lemat}
\label{sybb}  
Let $G=(V,\ka E)$ be a graph of siblings in a polybox code $V\subset S^d$, $u$ and $v$ be adjacent vertices  and let $d(u)=n$ and $d(v)=m$. If $n+m=2d$, then there are $i\in [d]$ and $l\in S$ such that   
\begin{equation}
\label{nmm3}
|V^{i,l}\cup V^{i,l'}|\geq 2d-2, 
\end{equation}
and if $n+m\leq 2d-1$, then
\begin{equation}
\label{nm3}
|V^{i,l}\cup V^{i,l'}|\geq n+m-1 
\end{equation}
for some $i\in [d]$ and $l\in S$. 
\end{lemat}

By $d(G)$ we denote the average degree of a graph $G$. 

\begin{lemat} 
\label{graph}
Let $G=(V,\ka E)$ be a simple graph, and let $m=\max\{d(v)+d(u):v,u\in V\;{\rm and}\; v,u\;{\rm are}\;{\rm adjacent}\}$. Then 
 $d(G)\leq m/2$.
\end{lemat}

Let $x\in ES$ and $i\in [d]$. Recall that 
$$
\pi^i_x=ES\times \cdots \times ES\times \{x\}\times ES\times \cdots \times ES,
$$
where $\{x\}$ stands at the $i$th position. If $V\subset (*S)^d$ is a polybox code, then the slice $\pi^i_x\cap \bigcup E(V)$ is a "flat" polybox in $(ES)^d$
(boxes which are contained in this polybox have the factor $\{x\}$ at the $i$th position (compare (\textbf{S}) in Section 3.1)). 
Therefore we define a polybox $(\pi^i_x\cap \bigcup E(V))_{i}$ in the $(d-1)$-box $(ES)^{d-1}$:
$$
(\pi^i_x\cap \bigcup E(V))_{i}=\bigcup \{(\breve{v})_{i}: v\in V\; {\rm and}\; \pi^i_x\cap \breve{v}\neq\emptyset\}.  
$$
The polybox $(\pi^i_x\cap \bigcup E(V))_{i}$ does not depend on a particular choice of a polybox code, because if $W$ is an equivalent polybox code to $V$, then $\bigcup E(V)=\bigcup E(W)$, and hence $(\pi^i_x\cap \bigcup E(V))_{i}=(\pi^i_x\cap \bigcup E(W))_{i}$. 


We will slice a polybox  $\bigcup E(V)$ by the set $\pi^i_x$ for various $x\in ES$. 
In particular, we will pay attention whether the polybox code $\{(v)_i:v\in V,\; \pi^i_x\cap \breve{v}\neq\emptyset\}$ is rigid or it contains a twin pair (compare \cite[Figure 6]{Kis}).

In \cite{Kis} we showed that any polybox code without twin pairs having at most seven words is rigid. Now we need a slightly better rigidity result: 

\begin{lemat}
\label{szt9}
If a polybox code $V\subset S^d$ does not contain a twin pair and $|V|\leq 9$, then it is rigid.
\end{lemat}
\proof Let $W$ be a polybox code equivalent to $V$ and suppose that $V\cap W=\emptyset$. It follows from Theorem \ref{12} that there is a twin pair $v,u$ in $W$. Let $v_i=u_i'$, and then $(v)_i=(u)_i$. We may assume that $(v)_i=b\cdots b$. In \cite[Lemma 3.11]{Kis} we showed that the lemma is true for $|V|\leq 7$. In the first part of the proof we show that the following implication holds: if $V$ is rigid for $|V|\leq 8$, then it is rigid for $|V|\leq 9$. So, assume that the lemma is true for $|V|\leq 8$.

Suppose first that $\breve{v}\cup \breve{u}\subset E(V^{i,a}\cup V^{i,a'})$. Since $(v)_i\sqsubseteq (V^{i,a})_i$ and  $(u)_i\sqsubseteq (V^{i,a'})_i$, by the assumption on $V$ and Statement \ref{stt} (a), we may assume that $|\{q\in V^{i,a}: \breve{q}\cap \breve{v}\neq\emptyset \}|=1$ and  $|\{q\in V^{i,a'}: \breve{q}\cap \breve{u}\neq\emptyset \}|\geq 5$. 
By ({\bf P}) in Section 3.1, there is $x\in \bigcup E(V^{i,a'})$ such that $(x)_i\not\in \bigcup E((V^{i,a})_i)$
and $q\in W$ such that $x\in \breve{q}$ and $q_i=a$. Suppose that $\breve{q}\cap \bigcup E(V^{i,b}\cup V^{i,b'})\neq\emptyset$. Since $q\sqsubseteq V$, this last set is an $i$-cylinder in the box $\breve{q}$ (compare Figure 7), and, by Statement \ref{stt} (c) and the fact that $|V|=9$, we have $|\{u\in V^{i,b}: \breve{u}\cap \breve{q}\neq\emptyset \}|\geq 2$ and $|\{u\in V^{i,b'}: \breve{u}\cap \breve{q}\neq\emptyset \}|\geq 2$. But then $|V|>9$, which is not true. Therefore, every such $x$ is contained in $\breve{q}$, where $q\in W$ and $q\sqsubseteq V^{i,a'}$. Thus, $|V^{i,a'}|\geq 6$ because if  $|V^{i,a'}|=5$, then, by Statement \ref{stt} (c),  the code $(V^{i,a'})_i$ covers only $(u)_i$. If  $|V^{i,a'}|=6$, then, by Statement \ref{stt} (b), $V^{i,a'}$ covers at most two words. Therefore, by (\textbf{V}) in Section 3.1,  $|V^{i,a}|\geq 4$, and thus  $|V|>9$, a contradiction. 

Let now $7\leq |V^{i,a'}|\leq 8$. Then we may assume that $|V^{i,b'}|=\emptyset$ and $|V^{i,b}|\leq 1$. 

Let $V^{i,b}=\{p\}$. Then $|V^{i,a'}|=7$ and $|V^{i,a}|=1$, and hence, by (\textbf{P}) in Section3.1, $W^{i,b}\neq\emptyset$. Since $u,v$ is a twin pair and $u_i=v_i'$, the polybox code $R=(W\setminus \{u,v\})\cup \{\bar{u},\bar{v}\}$, where $(\bar{u})_i=(u)_i$, $(\bar{v})_i=(v)_i$, $\bar{u}_i=\bar{v}_i'$ and $\bar{u}_i,\bar{v}_i\not\in \{a,a',b,b'\}$, is equivalent to $V$. Thus, by (\textbf{V}) in Section 3.1, $|R^{i,a'}|=6$ and $|R^{i,b}|=1$ and consequently $p\in W$, a contradiction.    

Let now $V^{i,b}=\emptyset$, and  let $x\in Ea'\cap Eb'$.  Then, by ({\bf S}) in Section 3.1, 
$$
\pi_x^i\cap \bigcup E(V)=\pi_x^i\cap \bigcup E(V^{i,a'}).
$$
Then $(\pi_x^i\cap \bigcup E(V^{i,a'}))_i=\bigcup_{w\in W}(\breve{w})_i$. By the inductive assumption,   $(\pi_x^i\cap \bigcup E(V^{i,a'}))_i$ is rigid, and thus, by (\textbf{P}) in Section 3.1, $w\in V$, for every $w\in W$ such that $x\in \breve{w}$, where $x\in \bigcup E(V^{i,a'})$ and $(x)_i\not \in \bigcup E((V^{i,a})_i)$. A contradiction. 

The implication: if $V$ is rigid for $|V|\leq 7$, then it is rigid for $|V|\leq 8$ is proved in the same manner.

We now consider the case $(\breve{v}\cup \breve{u})\subset \bigcup E(V^{i,a}\cup V^{i,a'} \cup V^{i,b}\cup V^{i,b'})$ and $(\breve{v}\cup \breve{u})\cap \bigcup E(V^{i,l}\cup V^{i,l'})\neq\emptyset$ for $l\in \{a,b\}$. The set $(\breve{v}\cup \breve{u})\cap E(V^{i,l}\cup V^{i,l'})$ is an $i$-cylinders in the $d$-box $\breve{v}\cup \breve{u}$ for $l\in \{a,b\}$. Therefore, $\bigcup \{(\breve{v})_i\cap (\breve{p})_i: p\in V^{i,l}\}=\bigcup \{(\breve{v})_i\cap (\breve{p})_i: p\in V^{i,l'}\}$ for $l\in \{a,b\}$ (recall that $(v)_i=(u)_i$ as $v$ and $u$ are twins). Let $F_l=\{(\breve{v})_i\cap (\breve{p})_i: p\in V^{i,l}\}$ for $l\in \{a,a',b,b'\}$.  By Statement \ref{stt} (c) and the fact that $|V|\leq 9$, we have $|F_l|\geq 2$ for every $l\in \{a,a',b,b'\}$. Suppose that $|F_l|=2$ for every $l\in \{a,a',b,b'\}$. Then, by Statement \ref{stt} (c), the sets $\bigcup F_a$ and $\bigcup F_b$ are $L$-shaped polyboxes (Figure 6a) and $\bigcup F_a\cup \bigcup F_b=(\breve{v})_i$. This is easy to see that a box (in this case $(\breve{v})_i$) cannot be represents as a sum of two disjoint $L$-shaped polyboxes. Therefore, since $|V|\leq 9$, we have $|F_a|=|F_{a'}|=2$ and    $|F_b|=2$, $|F_{b'}|=3$. 
By Statement \ref{stt} (c) there is a three-element set $A\subset [n]\setminus \{i\}$ such that:
$$
(V^{i,a})_{A^c}=\{aaa,ba'a\},\;(V^{i,a'})_{A^c}=\{a'a'a, aba\}
$$
$$
(V^{i,b})_{A^c}=\{a'aa,bba'\},\;\;(V^{i,b'})_{A^c}=\{a'ab,a'aa',ba'a'\}
$$ 
(compare Figure 6a and 6b).
Pick $x\in \bigcup E(V^{i,a})$ such that $(x)_i\not\in \bigcup E((V^{i,a'})_i)$. Then $x\in \breve{w}$, where $w\in W$ is, by (\textbf{P}) in Section 3.1, such that $w_i=a$, and consequently $\breve{w}\cap \breve{v}=\emptyset$ for $v\in V^{i,a'}$. Since $\breve{w}$ is dichotomous to the box $\breve{v}\cup \breve{u}$, there is $j\in [d]\setminus \{i\}$ such that $w'_j=v_j=u_j$. The sets $\{(\breve{v})_i\cap (\breve{q})_i: q\in V^{i,a}\cup V^{i,b}\}$ and $\{(\breve{v})_i\cap (\breve{q})_i: q\in V^{i,a}\cup V^{i,b'}\}$ are  dichotomous partition of $(\breve{v})_i$. Since $(v)_i=b\cdots b$, by Lemma \ref{34}, we have $(V^{i,a}\cup V^{i,b})_{A}=(V^{i,a}\cup V^{i,b'})_{A}=\{b\cdots b\}$. 

If $j\in A^c$, then $\breve{w}\cap \bigcup E(V^{i,a}\cup V^{i,a'} \cup V^{i,b}\cup V^{i,b'})=\emptyset$ (because $w_j=b'$) and then $|V|>9$ as $w\sqsubseteq V$. A contradiction. 

If $j\in A$, then $\breve{w}\cap \breve{q}=\emptyset$ for some $q\in V^{i,b}\cup V^{i,b'}$. This is easy to check that then the set $\breve{w}\cap \bigcup E( V^{i,b}\cup V^{i,b'})$ is not an $i$-cylinder in the box $\breve{w}$. On the other hand, this set has to be an $i$-cylinder because $w\sqsubseteq  V^{i,a}\cup V^{i,b}\cup V^{i,b'}$ and $w_i=a$ (compare Figure 7). 

Clearly, if $(\breve{v}\cup \breve{u})\subset \bigcup_{j\in [k]} E(V^{i,l_j}\cup V^{i,l_j'})$, where $k\geq 3$, then $|V|>9$.
\hfill{$\square$}

\medskip
In the next two lemmas we give  forbidden distributions of words in the considered codes $V$ and $W$.

\begin{lemat} 
\label{51}
If $V,W\subset S^d$, $S=\{a,a',b,b'\}$, are disjoint and equivalent polybox codes without twin pairs and $|V|=12$, then for every $i\in [d]$ the distribution 
$$
|V^{i,a}|=5,\;|V^{i,a'}|=1,\; |V^{i,b}|= 5,\; |V^{i,b'}|=1 
$$
is impossible.
\end{lemat}
\proof 
Assume on the contrary that there is $i\in [d]$ such that $|V^{i,a}|=5, |V^{i,a'}|=1$ and $|V^{i,b}|= 5, |V^{i,b'}|=1$. Let  $\{v\}=V^{i,a'}$ and $\{u\}=V^{i,b'}.$ If $u$ and $v$ is not an $i$-siblings, then the polybox code $\{(u)_i,(v)_i\}$ is rigid. Since $\pi_x\cap \bigcup E(V)=\pi_x\cap \bigcup E(W)$ for $x\in Ea'\cap Eb'$, it follows that $(v)_i=(w)_i$ and $(u)_i=(q)_i$, where $w\in W^{i,b'}$ and $q\in W^{i,a'}$. 
By (\textbf{Co}) in Section 3.1, we have $(w)_i\sqsubseteq (W^{i,b})_i$,  
$(q)_i\sqsubseteq (W^{i,a})_i$.
Thus, by Statement \ref{stt} (a), $|W^{i,b}|\geq 5$ and  $|W^{i,a}|\geq 5$. But $|W^{i,a'}|\geq 1$ and $|W^{i,b'}|\geq 1$, and therefore $|W^{i,b}|=5$.
Since $(w)_i\sqsubseteq (W^{i,b})_i$ and  $(w)_i\sqsubseteq (V^{i,a})_i$, the structure of $(W^{i,b})_i$ and $(V^{i,a})_i$ are such as given in  Statement \ref{stt} (a). Therefore, assuming without loss of generality that $(w)_i=b\cdots b$, we have  $(V^{i,a})_i=\{l_1l_2l_3,l_1'l_2'l_3',bl_2'l_3,l_1bl_3',l_1'l_2b\}$ and $(W^{i,b})_i=\{s_1s_2s_3,s_1's_2's_3',bs_2's_3,s_1bs_3',s_1's_2b\}$, where $l_i,s_i\in \{a,a'\}$ for $i=1,2,3$. It can be easily checked, using (\ref{dkostki}), that for every  $l_i,s_i\in \{a,a'\}$, $i=1,2,3$, there is  $y\in (ES)^d$ such that $(y)_i\in  \bigcup (E(V^{i,a}))_i\cap \bigcup (E(W^{i,b}))_i$ and $(y)_i\not\in (\breve{w})_i$. Observe now  that, again by (\ref{dkostki}), the point $y$  can be chosen such that $y_i\in Eb\setminus Ea$, and then $y\in  \bigcup E(W^{i,b})\setminus  \bigcup E(V^{i,a})$. Thus, $y\in \bigcup E(W)$ and $y\not\in \bigcup E(V)$, a contradiction. 

Before we consider the case when  $u,v$ is an $i$-siblings note that there is $w\in W$ such that $\breve{w}\cap \bigcup E(V^{i,a})\neq\emptyset$ and $\breve{w}\cap \bigcup E(V^{i,a'})\neq\emptyset$, for otherwise $V^{i,a}\sqsubseteq W^{i,a}$ and  $V^{i,a'}\sqsubseteq W^{i,a'}$. Then, by Statement \ref{stt} (b) and (a), respectively, $|W^{i,a}|\geq 7$ and $|W^{i,a'}|\geq 5$, and thus $|W^{i,a}|=7$ and $|W^{i,a'}|=5$, as $|V|=12$. Then $W^{i,b}\cup W^{i,b'}=\emptyset$, which is impossible, because, by (\textbf{P}) in Section 3.1, any $z\in \bigcup E(V^{i,b})$ such that $(z)_i\in \bigcup E((V^{i,b})_i)\setminus \bigcup E((V^{i,b'})_i)$ is covered by a box from $E(W^{i,b})$.  
Thus, the structures of $V^{i,a}$ and $V^{i,a'}$ are such as predicted in Statement \ref{stt} (c). In particular,  $(v)_i \sqsubseteq (V^{i,a})_i$.
 
Let now  $u$ and $v$ be an $i$-siblings. Then $(v)_i$ and $(u)_i$ are a twin pair. We can assume without loss of generality that $(v)_i=b\cdots b$. Since $(v)_i \sqsubseteq (V^{i,a})_i$, it follows that, by Statement \ref{stt} (a), $(V^{i,a})_{A^c}=\{l_1l_2l_3,l_1'l_2'l_3',bl_2'l_3,l_1bl_3',l_1'l_2b\}$,
where  $l_i\in \{a,a'\}$ for $i=1,2,3$, $A=\{i,i_1,i_2,i_3\}$ and $i_1,i_2,i_3$ are such as in (\ref{p5}) of Lemma \ref{34}. The words $(v)_i,(u)_i$ are a twin pair,  and therefore the polybox code $(V^{i,a})_i\cup (V^{i,b'})_i$ does not contain a twin pair, and hence, by Lemma \ref{szt9}, it is rigid. Then $(V^{i,a})_i\cup (V^{i,b'})_i=(W^{i,a})_i\cup (W^{i,b'})_i$, and consequently  $(V^{i,a})_i=(W^{i,b'})_i$ because $V$ and $W$ are disjoint. Thus, for $w\in W^{i,b'}$, by (\textbf{Co}) in Section 3.1, we have $(w)_i\sqsubseteq (V^{i,a'})_i$, and then $|V^{i,a'}|\geq 5$, a contradiction.
\hfill{$\square$}

\begin{lemat}
\label{rozk1}
Let $V,W\subset S^d$ be disjoint and equivalent polybox codes without twin pairs. If there are $i\in [d]$ and $l,s\in S$, $s\not\in \{l,l'\}$, such that $|V^{i,l}|=| V^{i,l'}|=1$ and  $|V^{i,s}|\neq |V^{i,s'}|$ or $|V^{i,l}|=1$ and $2\leq | V^{i,l'}|\leq 4$,  then $|V|>12$.
\end{lemat}
\proof 
By Theorem \ref{12}, $|V|\geq 12$. Suppose on the contrary that $|V|=12$. Let $|V^{i,l}|=|V^{i,l'}|=1$. By Statement \ref{stt} (d),  $V^{i,l}\sqsubseteq W^{i,l}$ and  $V^{i,l'}\sqsubseteq W^{i,l'}$, and thus, by Statement \ref{stt} (a),   $|W^{i,l}|\geq 5$ and  $|W^{i,l'}|\geq 5$.

Since $|V^{i,s}|\neq |V^{i,s'}|$, we may assume, by (\textbf{V}) in Section 3.1, that 
$W^{i,s}\neq\emptyset$. Consequently, we may assume that  $|W^{i,l}|=5$.

Suppose now that for every  $r\in S$, $r\not\in \{l,l'\}$ we have $W^{i,r}=\emptyset$ or  $W^{i,r'}=\emptyset$. Then, by (\ref{dkostki}), there is $x\in El$ such that $\bigcup E(W)\cap \pi^i_x= \bigcup E(W^{i,l})\cap \pi^i_x$. By Lemma \ref{szt9}, the polybox code $\{(w)_i: w\in W^{i,l}\}$ is rigid, and therefore $(v)_i=(w)_i$ for some $w\in W^{i,l}$ and $v\in V^{i,l}$ because   $\bigcup E(W)\cap \pi^i_x=\bigcup E(V)\cap \pi^i_x$. Thus, $v=w$, a contradiction. Hence, there is  $r\in S$, $r\not\in \{l,l'\}$ such that the sets  $W^{i,r}$ and $W^{i,r'}$ are non-empty. Clearly, $r=s$ and then $|W^{i,s}|=|W^{i,s'}|=1$ because $|W^{i,l}|\geq 5$, $|W^{i,l'}|\geq 5$ and $|V|=12$. It follows from Statement \ref{stt} (d) that $W^{i,s}\sqsubseteq V^{i,s}$ and  $W^{i,s'}\sqsubseteq V^{i,s'}$, and thus, by Statement \ref{stt} (a),   $|V^{i,s}|\geq 5$ and  $|V^{i,s'}|\geq 5$. Since , $|V^{i,s}|\neq |V^{i,s'}|$, we have $|V^{i,s}|+|V^{i,s'}|\geq 11$, and consequently $|V|>12$, a contradiction.

Let now  $|V^{i,l}|=1$ and $2\leq |V^{i,l'}|\leq 4$. By Statement \ref{stt} (d),  $V^{i,l}\sqsubseteq W^{i,l}$ and  $V^{i,l'}\sqsubseteq W^{i,l'}$, and by Statement \ref{stt} (a) and (b), respectively, $|W^{i,l}|\geq 5$ and $|W^{i,l'}|\geq  7$. Thus,  $|W^{i,l}|=5$ and $|W^{i,l'}|=7$ because $|V|=12$. This means that $W^{i,s}\cup W^{i,s'}=\emptyset$ for every $s\in S\setminus \{l,l'\}$, which implies, by (\textbf{C}) in Section 3.1, 
that the polybox codes  $(V^{i,s})_i$ and $(V^{i,s'})_i$ are equivalent. By Theorem \ref{12}, $|V|> |V^{i,s}|+|V^{i,s'}|\geq 24$, a contradiction.
\hfill{$\square$}

\medskip
The following will play lemma an important role in determining the structure of the polybox codes $V$ and $W$. 
\begin{lemat}
\label{ciecia}
Let $V,W\subset S^d$, where $S=\{a,a',b,b'\}$, be  disjoint and equivalent polybox codes without twin pairs. Assume that there are $i\in [d]$ and $l,p\in S, l\neq p,$ such that the sets  $V^{i,l}$ and $V^{i,p}$ are non-empty. If there is $x\in ES$ such that the polybox code $\{(v)_i: v\in V, \; \pi^i_x\cap \breve{v}\neq\emptyset\}$ does not contain a twin pair, then $|V|>12$.
\end{lemat}
\proof By Theorem \ref{12}, $|V|\geq 12$. For future applications, in the first part of the proof we do not assume that $S=\{a,a',b,b'\}$.

If $V^{i,s}=\emptyset$ or $V^{i,s'}=\emptyset$ for every $s\in S$, then, by (\textbf{P}) in Section 3.1, $\bigcup E(V^{i,l})=\bigcup E(W^{i,l})$, which implies that  $V^{i,l}$ and $W^{i,l}$ are equivalent. 
Similarly, $V^{i,p}$ and $W^{i,p}$ are equivalent. 
By Theorem \ref{12}, $|V|\geq |V^{i,l}|+|V^{i,p}|\geq 24$.

Let $p=l'$, and let $V= V^{i,l}\cup V^{i,l'}$. If $W= W^{i,l}\cup W^{i,l'}$, then the polybox codes $V^{i,l}$ and $W^{i,l}$ are equivalent and similarly, $V^{i,l'}$ and $W^{i,l'}$ are equivalent. Then, by Theorem \ref{12}, $|V|=|V^{i,l}|+|V^{i,l'}|\geq 24$.

If $W^{i,s}\neq\emptyset$ and $V= V^{i,l}\cup V^{i,l'}$, where $s\not\in \{l,l'\}$, then, by (\textbf{C}) in Section 3.1, 
the set $\bigcup E(W^{i,s}\cup W^{i,s'})$ is an $i$-cylinder, which gives $|W|\geq |W^{i,s}|+|W^{i,s'}|\geq 24$, and thus $|V|>12$.

Suppose now on the contrary that $|V|=12$. Furthermore,  we assume that $V^{i,l}\neq\emptyset$, $V^{i,l'}\neq\emptyset$ and $V^{i,s}=\emptyset$ or $V^{i,s'}=\emptyset$ for every $s\in S\setminus \{l,l'\}$.

Suppose that  $V^{i,s}\neq\emptyset$ for at least one $s\not\in \{l,l'\}$.
By (\ref{dkostki}) and (\textbf{S}) in Section 3.1, we can choose $x\in El\cap Es'$ and $y\in El'\cap Es'$ such that 
\begin{equation}
\label{111}
\pi^i_x\cap \bigcup E(V)= \pi^i_x\cap \bigcup E(V^{i,l})\; {\rm and} \;\pi^i_y\cap \bigcup E(V)= \pi^i_y\cap \bigcup E(V^{i,l'}).
\end{equation}
If $V^{i,l}$ is not rigid, then, by Lemma \ref{szt9}, $|V^{i,l}|\geq 10$. Since $|V|=12$, we have   $|V^{i,l}|=10$,  $|V^{i,l'}|=1$ and  $|V^{i,s}|=1$. Then $V^{i,s}\sqsubseteq W^{i,s}$, and consequently, by Statement \ref{stt} (a), $|W^{i,s}|\geq 5$. Therefore, for every $x\in El\cap Es'$ and $y\in El'\cap Es$ we have 
$$
 |(\pi^i_{x}\cap \bigcup E(V))_i|_0\geq 10\; {\rm and}\; |(\pi^i_y\cap \bigcup E(V))_i|_0\geq 5,     
$$
where $|F|_0$ is the number of boxes in arbitrary proper suit for $F$ (see Section 2).
Now it is easy to see (compare \cite[Lemma 3.6]{Kis}) that $|V|\geq 15$, which contradicts the assumption that $|V|=12$.

If the codes $V^{i,l}$ and $V^{i,l'}$ are rigid, then $W^{i,l}\cup W^{i,l'}=\emptyset$, for otherwise taking $w\in W^{i,l}$ we get, by (\ref{111}) and the rigidity of $V^{i,l}$, $(v)_i=(w)_i$ for some $v\in V^{i,l}$, and thus $v=w$, a contradiction. Since $W^{i,l}\cup W^{i,l'}=\emptyset$, the set $\bigcup E(V^{i,l}\cup V^{i,l'})$ is, by (\textbf{C}) in Section 3.1, an $i$-cylinder, and since $V^{i,l}$ and $V^{i,l}$ are rigid, the set $V^{i,l}\cup V^{i,l'}$ consists of twin pairs, which is a contradiction.

Thus, we may assume that the sets $V^{i,l}, V^{i,l'}, V^{i,s}$ and $V^{i,s'}$ are non-empty.

Let now $S=\{a,a',b,b'\}$. (We still assume that $|V|=12$.) It follows from the above that $V^{i,l}\neq\emptyset$ for $l\in \{a,a',b,b'\}$.  

Let $x\in Ea\cap Eb$ be such that the set $\{(v)_i: v\in V,\; \pi^i_x\cap \breve{v}\neq\emptyset\}$ does not contain a twin pair.  

If $|(\pi^i_x\cap \bigcup E(V))_i|_0\geq 10$, then $|(\pi^i_x\cap \bigcup E(V))_i|_0=10$ and  $|(\pi^i_y\cap \bigcup E(V))_i|_0=2$ for $y\in Ea'\cap Eb'$, and consequently $|V^{i,a'}|=|V^{i,b'}|=1$. Since $V=V^{i,a}\cup V^{i,a'}\cup V^{i,b}\cup V^{i,b'}$, by Lemma \ref{rozk1},  $|V^{i,a}|=5$, $|V^{i,a'}|=1$ and  $|V^{i,b}|=5$, $|V^{i,b'}|=1$, which is, by Lemma \ref{51}, impossible.

If $|(\pi^i_x\cap \bigcup E(V))_i|_0\leq 9$, then, by Lemma \ref{szt9}, the polybox code $\{(v)_i: v\in V,\; \pi^i_x\cap \breve{v}\neq\emptyset\}$ is rigid. Thus, there are  $v\in V^{i,a}, u\in V^{i,b}$ and $w,p\in W$ with $w_i\in \{b,b'\}$ and $p_i\in \{a,a'\}$ such that $(v)_i=(w)_i$ and $(u)_i=(p)_i$. Then, by (\textbf{Co}) in Section 3.1 and Statement \ref{stt} (a), $|V^{i,a'}|\geq 5$, $|V^{i,b'}|\geq 5$. Since $|V|=12$, we have $|V^{i,a'}|=5$, $|V^{i,a}|=1$ and  $|V^{i,b'}|=5$, $|V^{i,b}|=1$, which is, by Lemma \ref{51}, impossible. Thus $|V|>12$.
\hfill{$\square$}

\medskip
We now once again indicate a forbidden distribution of words in $V$ and $W$.   

\begin{lemat} 
\label{33}
If $V,W\subset S^d$, where $d\geq 5$ and $S=\{a,a',b,b'\}$, are disjoint and equivalent polybox codes without twin pairs and $|V|=12$, then the distributions of words in $V$ of the forms
$$
|V^{i,a}|=|V^{i,a'}|=|V^{i,b}|=|V^{i,b'}|=3\;\; {\rm for\;\; every}\;\; i\in [d]
$$
is impossible. 
\end{lemat}
\proof By Lemma \ref{ciecia}, we may assume that for every $i\in [d]$ and $l,s\in \{a,a',b,b'\}, l\not\in \{s,s'\}$, there are $i$-siblings in the set $V^{i,l}\cup V^{i,s}$.


Suppose on the contrary that $V$ has the distribution  $|V^{i,l}|=|V^{i,l'}|=3$  for every $i\in [d]$ and $l\in \{a,b\}$. Let $G=(V,\ka E)$ be a graph of siblings in $V$. Note that, it follows from the assumption on $i$-siblings in $V$ that $|\ka E|\geq 4d$.

Let $u^0,v^0\in V$ be such that 
$$
d(v^0)+d(u^0)=\max\{d(v)+d(u):v,u\in V\;{\rm and}\; v,u\;{\rm are}\;{\rm adjacent}\}.     
$$  
If $d(v^0)+d(u^0)\geq 8$, then, by Lemma \ref{sybb}, there are $i\in [d]$ and $l\in \{a,b\}$ such that $|V^{i,l}\cup V^{i,l'}|\geq 7$, which contradicts the assumption on the distribution of words in $V$. 
On the other hand, if $d(v^0)+d(u^0)\leq 6$, then, by Lemma \ref{graph}, $d(G)\leq 3$. Since $|\ka E|\geq 20$, we have $|V|>12$, which is a contradiction. 

Let $d(v^0)+d(u^0)=7$. It can be easily shown that there are $i,j\in [d]$, $i\neq j$, and $l,s\in \{a,b\}$ such that $|V^{i,l}\cup V^{i,l'}|=6$ and  $|V^{j,s}\cup V^{j,s'}|=6$, where $V^{i,l},V^{i,l'}, V^{j,s}, V^{j,s'}\subset N(u^0)\cup N(v^0)$ (compare \cite[Lemma 3.9]{Kis}). We can assume without loss generality, as $i\neq j$, that $l=s=a$. Let $\{u\}=(N(u^0)\cup N(v^0))\setminus  (V^{i,a}\cup V^{i,a'})$ and $\{v\}=(N(u^0)\cup N(v^0))\setminus  (V^{j,a}\cup V^{j,a'})$. Clearly, $u_i\in \{b,b'\}$ and $v_j\in \{b,b'\}$.  Thus, $u\neq v$, for otherwise $u\not\in N(u^0)\cup N(v^0)$. Assume without loss of generality that $u_i=b$ and $v_i=a$. Since $w_i,w_j\in \{b,b'\}$ for every $w\in V\setminus (N(u^0)\cup N(v^0))$ and $p_i,p_j\in \{a,a'\}$ for every $p\in (N(u^0)\cup N(v^0))\setminus \{u,v\}$, a vertex $w\in  V\setminus (N(u^0)\cup N(v^0))$ can be joined only with $u$ or $v$. This means that there is no $i$-siblings $q,t\in V$ such that $q_i=b'$ and $t_i=a'$, a contradiction. 

\hfill{$\square$}

\medskip
In the next lemma we show that the polyboxes $V$ and $W$ can be written down in the alphabet $S=\{a,a',b,b'\}$.

\begin{lemat}
\label{abc}
Let $V,W\subset S^d$ be disjoint and equivalent polybox codes without twin pairs, and let $V$ be extensible to a partition code. If there is $i\in [d]$ such that $V^{i,l}\cup V^{i,l'}\neq\emptyset$ for at least three $l\in S$, then $|V|> 12$.
Thus, if $|V|=|W|=12$, then $V,W\subset \{a,a',b,b'\}^d$.  
\end{lemat}
\proof It follows from Theorem \ref{12} that $|V|\geq 12$. Suppose on the contrary that  $|V|=12$. 

In the same way as in the first part of the proof of Lemma \ref{ciecia} we show that for every $i\in [d]$  there are at least two letters $l,s\in S, $ $l\not\in \{s,s'\}$, such that 
\begin{equation}
\label{2}
V^{i,l}, V^{i,l'},V^{i,s}, V^{i,s'}\neq\emptyset  
\end{equation}

If $W^{i,l}\cup W^{i,l'}=\emptyset$, then, by (\textbf{C}) in Section 3.1 and Theorem \ref{12}, we have 
 $|V|>|V^{i,l}|\cup |V^{i,l'}|\geq 24$, a contradiction. Thus, $W^{i,l}\cup W^{i,l'}\neq \emptyset$ and $W^{i,s}\cup W^{i,s'}\neq\emptyset$

Suppose that $V^{i,r}\neq\emptyset$ and $V^{i,r'}=\emptyset$ for some $r\in S\setminus \{l,l',s,s'\}$. By (\textbf{P}) in Section 3.1, for every $w\in W$ such that $\breve{w}\cap \bigcup E(V^{i,r})\neq\emptyset$ we have $w_i=r$, and hence $V^{i,r}\sqsubseteq W^{i,r}$ from where, by Statement \ref{stt} (a), we obtain $|W^{i,r}|\geq 5$. 
Then $|W^{i,l}\cup W^{i,l'}|\leq 3$ or $|W^{i,s}\cup W^{i,s'}|\leq 3$ because $|W|=12$. 

Let $|W^{i,l}\cup W^{i,l'}|\leq 3$ and  $W^{i,l}\neq\emptyset$, $W^{i,l'}\neq\emptyset$. Since  $|W^{i,r}|\geq 5$ and $W^{i,s}\cup W^{i,s'}\neq\emptyset$, by Lemma \ref{rozk1}, $|W|>12$, a contradiction. 
    
If  $|W^{i,l}\cup W^{i,l'}|\leq 3$ and $W^{i,l}\neq\emptyset,W^{i,l'}=\emptyset$, then  $W^{i,l}\sqsubseteq V^{i,l}$, and by  Statement \ref{stt} (a), $|V^{i,l}|\geq 5$. Moreover, since $W^{i,l'}=\emptyset$, by (\textbf{P}) in Section 3.1, we have  
$$
(V^{i,l'})_i\sqsubseteq (V^{i,l})_i. 
$$
Let $|V^{i,l'}|=1$ and $|V^{i,l}|\geq 5$. Then, by (\textbf{V}) in Section 3.1,  $|W^{i,l}|\geq 4$.
Similarly, if $|V^{i,l'}|\in \{2,3\}$, then,   by Statement \ref{stt} (b), $|V^{i,l}|\geq 7$, and consequently $|W^{i,l}|\geq 4$. In both cases we get a contradiction to $|W^{i,l}\cup W^{i,l'}|\leq 3$.

If $|V^{i,l'}|\geq 4$ and $|V^{i,l}|\geq 7$, then, by (\ref{2}), $|V|>12$, a contradiction.

Therefore, $V^{i,r}, V^{i,r'}\neq\emptyset$. Moreover, we showed that for every $j\in [d]$ and $p\in S$
\begin{equation}
\label{1}
V^{j,p}\neq\emptyset \Leftrightarrow V^{j,p'}\neq\emptyset.
\end{equation}

Suppose that  $|V^{i,l}|=|V^{i,l'}|=1$. It follows from Statement \ref{stt} (c) that $V^{i,l}\sqsubseteq W^{i,l}$ and  $V^{i,l'}\sqsubseteq W^{i,l'}$. By Statement \ref{stt} (a), we have $|W^{i,l}|\geq 5$ and $|W^{i,l'}|\geq 5$. Clearly, similarly like in the distribution (\ref{2}), there are $s_1,r_1\in S\setminus\{l,l'\}$, $s_1\not\in \{r_1,r_1'\}$, such that the sets $W^{i,s_1}, W^{i,s_1'}, W^{i,r_1}$ and $W^{i,r_1'}$ are non-empty. Thus, $|W|>12$, a contradiction. 

If  $|V^{i,l}|=1$ and $|V^{i,l}|=2$, then, by Lemma \ref{rozk1}, $|V|>12$, a contradiction. 

Thus $|V^{i,p}\cup V^{i,p'}|=4$, and by Lemma \ref{rozk1}, $|V^{i,p}|=|V^{i,p'}|=2$ for $p\in \{l,s,r\}$.

Let $d=4$. Since $V$ can be extended to a partition code $U\subset S^4$, i.e. $|U|=16$, we have 
$V^{i,p}\cup V^{i,p'}\subset U^{i,p}\cup U^{i,p'}$ for $p\in \{l,s,r\}$. As $|U^{i,p}|=|U^{i,p'}|$ and $\bigcup (E(V^{i,p}))_i\neq \bigcup (E(V^{i,p'}))_i$  for $p\in \{l,s,r\}$, at least two words are needed to complete the set   $V^{i,p}\cup V^{i,p'}$ to the set $U^{i,p}\cup U^{i,p'}$ for $p\in \{l,s,r\}$. But then $|U^{i,p}\cup U^{i,p'}|\geq 6$ for $p\in \{l,s,r\}$, and thus $|U|>16$ which is a contradiction.

Let now $d\geq 5$. For every $p\in \{l,s,r\}$ there is $w\in W$ such that  $\breve{w}\cap \bigcup E(V^{i,p})\neq\emptyset$ and $\breve{w}\cap \bigcup E(V^{i,p'})\neq\emptyset$, for otherwise  $V^{i,p}\sqsubseteq W^{i,p}$ and $V^{i,p'}\sqsubseteq W^{i,p'}$ for some $p$, and thus, by Statement \ref{stt} (b), $|W^{i,p}|\geq 7$ and $|W^{i,p'}|\geq 7$. Then $|W|>12$, a contradiction.

By Statement \ref{stt} (c), for every $p\in \{l,s,r\}$ there are $i_1(p),i_2(p)\in [d]\setminus \{i\}, i_1(p)<i_2(p)$, such that 
$$
(V^{i,p})_{A^c}=\{l_1(p)'l_2(p),l_1(p)s_2(p)\},\; \;\; (V^{i,p'})_{A^c}=\{l_1(p)l_2(p)',s_1(p)l_2(p)\}
$$ 
and 
$$
(V^{i,p})_{A\cup \{i\}}=(V^{i,p'})_{A\cup \{i\}}=\{o(p)\}
$$
where $o(p)\in S^{d-3}$ for $p\in \{l,s,r\}$, $A=\{i_1(p),i_2(p)\}$ and $l_i(p),s_i(p)\in S$, $l_i(p)\not \in \{s_i(p),s_i(p)'\}$ for $i=1,2$.
 
Let $j\in [d]\setminus \{i\}$ be such that  $v_j=o_j(p)$ for every $v\in V^{i,p}\cup V^{i,p'}$ and $p\in \{l,s,r\}$. Then  $V^{j,o_j(p)}\neq\emptyset$ and  $V^{j,o_j(p)'}=\emptyset$ for some $p\in \{l,s,r\}$, which contradicts (\ref{1}) or $V= V^{j,o_j(p)}\cup V^{j,o_j(p)'}$, which contradicts (\ref{2}). 

Since $d\geq 5$, there is $j\in [d]\setminus \{i\}$ and there are two letters in $S$, say $l$ and $s$, such that
$$
v_j=o_j(l)\;\; {\rm for}\;\; v\in V^{i,l}\cup V^{i,l'}\;\; {\rm and}\;\; v_j=o_j(s)\;\; {\rm for}\;\; v\in V^{i,s}\cup V^{i,s'},
$$
where $o_j(l),o_j(s)\in S$ and 
$$
v_j\in \{l_1(r), l_1(r)', s_1(r)\}\;\; {\rm for}\;\; v\in V^{i,r}\cup V^{i,r'},
$$
or
$$
v_j\in \{l_2(r), l_2(r)', s_2(r)\}\;\; {\rm for}\;\; v\in V^{i,r}\cup V^{i,r'},
$$
where   $l_k(r)\not \in \{s_k(r),s_k(r)'\}$ for $k=1,2$. We consider the first case (the second case is considered in the same manner).

Let $o_j(l)=o_j(s)$. Then $|V^{j,o_j(l)}|=8$. By (\ref{1}), we have $V^{j,o_j(l)'}\neq\emptyset$, and by (\ref{2}) we have $V^{j,s},V^{j,s'}\neq\emptyset$ for at least one $s\not\in \{o(j),o(j)'\}$. Thus, by Lemma \ref{rozk1}, $|V|> 12$, a contradiction. 

Let $o_j(l)=o_j(s)'$.  If $l_1(r)\in \{o_j(l), o_j(l)'\}$, then $V^{j,s_1(r)'}=\emptyset$, which contradicts (\ref{1}), as  $V^{j,s_1(r)}\neq\emptyset$. If  $l_1(r)\not\in \{o_j(l), o_j(l)'\}$, then $|V^{j,l_1(r)}|=2$, $|V^{j,l_1(r)'}|=1$ and, by Lemma \ref{rozk1}, $|V|> 12$, a contradiction.

Finally, let $o_j(l)\not\in \{o_j(s), o_j(s)'\}$. 
If $s_1(r)\not\in \{o_j(l), o_j(l)',o_j(s), o_j(s)'\}$, then  $V^{j,s_1(r)'}=\emptyset$, which contradicts (\ref{1}). Let $s_1(r)\in \{o_j(l), o_j(l)'\}$. Then  $l_1(r)\not\in \{o_j(l), o_j(l)'\}$, and thus $V^{j,s_1(r)'}=\emptyset$ or $|V^{j,s_1(r)}|=1$ and $|V^{j,s_1(r)'}|=4$. The first case contradicts (\ref{1}), and in the second case, by Lemma \ref{rozk1}, $|V|>12$, which is also a contradiction.  

To show that  $V,W\subset \{a,a',b,b'\}^d$ assume on the contrary that  $V\subset \{a,a',b,b'\}^d$ and  $W\subset \{c,c',d,d'\}^d$, where $\{a,a',b,b'\}\neq \{c,c',d,d'\}$. Let $V^{i,c}\cup V^{i,c'}=\emptyset$. By (\textbf{C}) in Section 3.1, the set  $\bigcup E(W^{i,c}\cup W^{i,c'})$ is an $i$-cylinder and consequently $|W^{i,c}|+|W^{i,c'}|\geq 24$, a contradiction.
\hfill{$\square$}

\medskip
At the end of this section we show that the computations can be  made mainly for $d=4,5$ and only in one case for $d=6$.

\begin{lemat}
\label{6}
Let $V,W\subset S^6$, where $S=\{a,a',b,b'\}$, be disjoint and equivalent polybox codes without twin pairs such that the distribution of words in $V$ is different from $|V^{i,a}|=|V^{i,a'}|=5$ and $|V^{i,b}|=|V^{i,b'}|=1$ for every $i\in [6]$. Then $|V|> 12$. 
\end{lemat}
\proof By Theorem \ref{12}, $|V|\geq 12$. Let $G=(V,\ka E)$ be a graph of siblings in $V$. 
By Lemma \ref{ciecia}, we assume that for every $i\in [d]$ and $l,s\in \{a,a',b,b'\}, l\not\in \{s,s'\}$, there is an $i$-siblings in the set $V^{i,l}\cup V^{i,s}$. Thus, for every $i\in [6]$ and every $\{l,s\}\in \{\{b,a\},\{b,a'\},\{b',a\},\{b',a'\}\}$ there is an edge $(v,u)\in \ka E$ such that $\{v_i,u_i\}=\{l,s\}$. In particular, for every $i\in [6]$ there are at least $4$ edges with the colour $i$, and therefore $|\ka E|\geq 24$.  

Let
$u^0,v^0\in V$ be such that 
$$
d(v^0)+d(u^0)=\max\{d(v)+d(u):v,u\in V\;{\rm and}\; v,u\;{\rm are}\;{\rm adjacent}\}.   
$$ 
We may assume without loss of generality that  $u^0=aaaaaa$, $v^0=ba'aaaa$.

By the assumption on the distribution of words in $V$ and Lemma \ref{rozk1},
\begin{equation}
\label{roz}
|V^{i,l}\cup V^{i,l'}|\geq 4
\end{equation}
for every $i\in [6]$ and $l\in \{a,b\}$.
 
If $d(u^0)+d(v^0)\geq 10$, then, by Lemma \ref{sybb}, there are $i\in [6]$ and $l\in \{a,b\}$ such that $|V^{i,l}\cup V^{i,l'}|\geq 9$. By (\ref{roz}), $|V|>12$. 

If  $d(u^0)+d(v^0)=9$, then it can be easily seen that there are $i,j\in [6]$, $i\neq j$, and $l,s\in \{a,b\}$ such that $|V^{i,l}\cup V^{i,l'}|\geq 8$ and $|V^{j,s}\cup V^{j,s'}|\geq 8$ (compare \cite[Lemma 3.9]{Kis}). Hence, by (\ref{roz}), $|V^{i,l}\cup V^{i,l'}|=8$ and $|V^{j,s}\cup V^{j,s'}|=8$. If $|V|=12$, then along the same lines as in  the proof of the second part of Lemma \ref{51} (the case  $d(u^0)+d(v^0)=7$) we show that there is $i\in [6]$ such that the set $\ka E$ contains less than $4$ edges with the colour $i$, which contradicts the assumption on $\ka E$.

Suppose on the contrary that $|V|=12$. 

Let $d(u^0)+d(v^0)=8$, $\{w^1,w^2,w^3,w^4\}=V\setminus (N(u^0)\cup N(v^0))$ and
$$
 |(V^{i,b}\cup V^{i,b'})\cap  (N(u^0)\cup N(v^0))|=n_i.
$$
Assume first that there is $i\in \{2,\ldots ,6\}$, say $i=6$, such that $n_6=0$. Then, by (\ref{roz}), $\{w^1,w^2,w^3,w^4\}=V^{6,b}\cup V^{6,b'}$. There are at least 4 edges $(w,v)$ with the colour $6$ in $\ka E$, and it follows from the assumption $n_6=0$ that it must be $w\in V^{6,b}\cup V^{6,b'}$ and $v\in N(u^0)\cup N(v^0)$.

Note that, there is $w\in W$ such that $\breve{w}\cap \bigcup E(V^{6,b})\neq\emptyset$ and $\breve{w}\cap \bigcup E(V^{6,b'})\neq\emptyset$, for otherwise $V^{6,b}\sqsubseteq W^{6,b}$ and  $V^{6,b'}\sqsubseteq W^{6,b'}$, and consequently, by Statement \ref{stt} $(a)$ and $(b)$, $|W^{6,b}|+|W^{6,b'}|\geq 12$. Since $W^{6,a}\cup W^{6,a'}\neq\emptyset$, we have $|W|>12$, which is a contradiction. By Lemma \ref{rozk1}, $|V^{6,b}|=|V^{6,b'}|=2$, and by Statement \ref{stt} $(c)$, 
$(V^{6,b})_{A}=(V^{6,b'})_{A}=\{l_4l_5\}$, $l_4,l_5\in S$, where we assumed without loss of generality that $A^c=\{4,5\}$. 

Let $w^i_4,w^i_5\in \{b,b'\}$ for $i\in \{1,2,3,4\}$. 
Since $n_6=0$, an edge $(w,v)$ has the colour $6$ and then $v_4,v_5\in \{b,b'\}$, and consequently $v\not\in N(u^0)\cup N(v^0)$. Thus, there are no edges of the colour 6 in $\ka E$, a contradiction.

Let now $w^i_5\in \{a,a'\}$ for $i\in \{1,2,3,4\}$. Then, by (\ref{roz}),  $n_5\geq 4$. Note that, if $v_i,v_j\in \{b,b'\}$, $i,j\in \{2,3,4,5\}, i\neq j$, then $v\not\in N(u^0)\cup N(v^0)$. Thus, $n_i=0$ for some $i\in \{2,3,4\}$, and consequently an edge $(w,v)$, where $w\in \{w^1,w^2,w^3,w^4\}$ and $v\in N(u^0)\cup N(v^0)$, has the colour $i$. On the other hand $(w,v)$ must be of the colour $6$ as $n_6=0$, a contradiction.   

Thus, $n_i\geq 1$ for every $i\in [6]$, and therefore it suffices to consider two cases: $n_2=2$, $n_3=\cdots =n_6=1$ and $n_1=\cdots =n_6=1$. 

It follows from (\ref{roz}) that in the first case there are at least two words in the set  $\{w^1,w^2,w^3,w^4\}$, say these are $w^1$ and $w^2$, which have the letters $b$ or $b'$ at at least three positions $i,j,k\in \{3,4,5,6\}$, which means that they cannot be adjacent to vertices from the set $N(u^0)\cup N(v^0)$. In the second case there are at least three such words; assume that these are $w^1,w^2$ and $w^3$. It is easy to verify that in the first case there are at most $8$ edges with ends in the set $N(u^0)\cup N(v^0)$, and in the second case there are at most $12$ such edges. The maximal number of edges with ends in $\{w^1,w^2,w^3,w^4\}$ is four as the graph $G$ does not contain triangles. Thus, in the second case in order to obtain $|\ka E|\geq 24$ it must be $d(w^4)\geq 8$, which is impossible since $d(v)\leq 6$ for every $v\in V$. For the same reason in the first case it must be $d(w^3)=d(w^4)=6$ and the vertices $w^3$, $w^4$ are adjacent to vertices $v\in (N(u^0)\cup N(v^0))\setminus \{u^0,v^0\}$. Thus, the vertices $w^3$ and $w^4$ are not adjacent to $w^1$ and $w^2$ and then $|\ka E|<24$, which contradicts the assumption on $\ka E$. 

Let now $d(v^0)+d(u^0)\leq 7$.  It follows from Lemma \ref{graph} that $d(G)\leq 7/2$, 
and since $d(G)|V|=2|\ka E|$ and $2|\ka E|\geq 48$, we have $|V|>12$.
\hfill{$\square$}

\section{Computations}
In this section we describe the computations which lead to the determination of all possible equivalent and disjoint polybox codes $V,W\subset S^d$ without twin pairs, having 12 words each, where $S=\{a,a',b,b'\}$ and $d\in \{4,5,6\}$. The structure of such polybox codes $V,W$ is given in Theorem \ref{s12}. 
   
The longest part of the paper is devoted to the preparations of computations. It seems hopeless to make the computations without any initial configurations of words, where by an initial configuration of words we mean a some number of words or their fragments in the constructing code $V$ (see tables in this section).

An immediate consequence of lemmas \ref{51}--\ref{abc} is the following result on the distribution of words in $V$:

\begin{wn}
\label{rozk}
Let $V,W\subset S^d$, $S=\{a,a',b,b'\},$ be disjoint and equivalent polybox codes without twin pairs, and let $|V|=12$. Then for every $i \in [d]$  the distribution of words in $V$ takes one of the forms:

$$
1.\;\; |V^{i,a}|=7,\;|V^{i,a'}|=1,\; |V^{i,b}|= 2,\; |V^{i,b'}|=2
$$
$$
2.\;\; |V^{i,a}|=6,\;|V^{i,a'}|=2,\; |V^{i,b}|= 2,\; |V^{i,b'}|=2
$$
$$
3.\;\; |V^{i,a}|=6,\;|V^{i,a'}|=1,\; |V^{i,b}|= 2,\; |V^{i,b'}|=3
$$
$$
4.\;\; |V^{i,a}|=5,\;|V^{i,a'}|=3,\; |V^{i,b}|= 2,\; |V^{i,b'}|=2
$$
$$
5.\;\; |V^{i,a}|=5,\;|V^{i,a'}|=2,\; |V^{i,b}|= 2,\; |V^{i,b'}|=3
$$
$$
6.\;\; |V^{i,a}|=5,\;|V^{i,a'}|=1,\; |V^{i,b}|= 3,\; |V^{i,b'}|=3
$$
$$
7.\;\; |V^{i,a}|=5,\;|V^{i,a'}|=1,\; |V^{i,b}|= 2,\; |V^{i,b'}|=4
$$
$$
8.\;\; |V^{i,a}|=5,\;|V^{i,a'}|=5,\; |V^{i,b}|= 1,\; |V^{i,b'}|=1
$$
$$
9.\;\; |V^{i,a}|=4,\;|V^{i,a'}|=4,\; |V^{i,b}|= 2,\; |V^{i,b'}|=2
$$
$$
10.\;\; |V^{i,a}|=4,\;|V^{i,a'}|=3,\; |V^{i,b}|= 2,\; |V^{i,b'}|=3
$$
$$
11.\;\; |V^{i,a}|=4,\;|V^{i,a'}|=2,\; |V^{i,b}|= 2,\; |V^{i,b'}|=4
$$
$$
12.\;\; |V^{i,a}|=3,\;|V^{i,a'}|=3,\; |V^{i,b}|= 2,\; |V^{i,b'}|=4.
$$
If $d=4$, then
$$
13.\;\; |V^{i,a}|=3,\;|V^{i,a'}|=3,\; |V^{i,b}|= 3,\; |V^{i,b'}|=3.
$$

Moreover, in every case $1-13$, except for the case $8$, for every $l\in \{a,b\}$ there is $w\in W$ such that $\breve{w}\cap \bigcup E(V^{i,l})\neq\emptyset$ and  $\breve{w}\cap \bigcup E(V^{i,l'})\neq\emptyset$. 

\end{wn}
\proof We prove the second part of the corollary. For every $w\in W$ we have $w\sqsubseteq V$. It follows from Lemma \ref{=c} that the family of pairwise dichotomous boxes $\{\breve{w}\cap \breve{v}\neq\emptyset: v\in V\}$ is a suit for the $d$-box $\breve{w}$ without twin pairs.  Suppose on the contrary that for every $w\in W$ at most one of the sets $\breve{w}\cap \bigcup E(V^{i,l})$ and  $\breve{w}\cap \bigcup E(V^{i,l'})$ is non-empty for some $l\in \{a,b\}$. Then $V^{i,l}\sqsubseteq W^{i,l}$ and $V^{i,l'}\sqsubseteq W^{i,l'}$. Since $|V^{i,l}|\geq 2$ or  $|V^{i,l'}|\geq 2$ (recall that the case $8$ has been excluded), we have, by Statement \ref{stt} (b), $|W^{i,l}|\geq 7$ or  $|W^{i,l'}|\geq 7$. But $|W^{i,l}|\geq 5$ and  $|W^{i,l'}|\geq 5$, by Statement \ref{stt} (a). Thus, $|W^{i,l}\cup W^{i,l'}|\geq 12$, and as  $|W^{i,s}\cup W^{i,s'}|> 0$, $s\not\in \{l,l'\}$,  we have $|W|>12$, a contradiction.
\hfill{$\square$}

\subsection{Initial configurations of the words.} It follows from Corollary \ref{rozk} that there is $w\in W$ such that $\breve{w}\cap \bigcup E(V^{i,b})\neq\emptyset$ and  $\breve{w}\cap \bigcup E(V^{i,b'})\neq\emptyset$ when the polybox code $V$ has the distributions of words $1-5$ and $9-13$ and there is $w\in W$ such that $\breve{w}\cap \bigcup E(V^{i,a})\neq\emptyset$ and  $\breve{w}\cap \bigcup E(V^{i,a'})\neq\emptyset$ when $V$ has the distributions $6$ and $7$.

We now show how the initial configurations of words are established. It is done in detail for the distributions $1,2,4$ and $9$.
The remaining configurations are determined in a similar manner. 

Let $V$ has the distribution of words of the form $1,2,4$ or $9$ of Corollary \ref{rozk}. Then, by Statement \ref{stt} (d), there is $w\in W$ such that  $|\{\breve{w}\cap \breve{v}\neq\emptyset: v\in V^{i,b}\}|=2$ and $|\{\breve{w}\cap \breve{v}\neq\emptyset: v\in V^{i,b'}\}|=2$. By Statement \ref{stt} (c), the set $\breve{w}\cap \bigcup E(V^{i,b})\cup \breve{w}\cap \bigcup E(V^{i,b'})$ is an $i$-cylinder in the $d$-box $\breve{w}$. 
Again by Statement \ref{stt} (c), there are $i_1,i_2\in [d]\setminus \{i\},i_1<i_2$, such that 
$$
(V^{i,b})_{A^c}=\{w_{i_1}l_2,l_1l_2'\},\;\; (V^{i,b'})_{A^c}=\{l_1l_2,l_1'w_{i_2}\}
$$
and 
$$
(V^{i,b})_{\{i,i_1,i_2\}}=(V^{i,b'})_{\{i,i_1,i_2\}}=\{p\},
$$
where $A=\{i_1,i_2\}$, $l_j\not\in \{w_{i_j},w_{i_j}'\}$ for $j=1,2$, and $p\in S^{d-3}$ (compare Figure 7).  Clearly, without loss of generality we may assume that $w_{i_1}=w_{i_2}=b$ and $p=b\cdots b$ (recall that $V,W\subset \{a,a',b,b'\}^d$).  Then, by Statement \ref{stt} (c), $l_1,l_2\in \{a,a'\}$, so we may assume that $l_1=l_2=a$.

Let $V$ be such as in the case $1$ of Corollary \ref{rozk}, and let $i=1,i_1=2,i_2=3$ and $d=5$. For the computations we take $a=+1,a'=-1,b=+2$ and $b'=-2$. We arrange the words from the set $V^{1,+2}\cup V^{1,-2}$ as the first four rows of the matrix $M$ of size $12\times 5$. The remaining eight words from  $V^{1,+1}\cup V^{1,-1}$ are unknown; we only know their first letters, i.e., $+1$ and $-1$ and that $|V^{1,+1}|=7$ and $|V^{1,-1}|=1$. The matrix $M$ has a form

\begin{displaymath}
{\rm M}=
\left[\begin{array}{ccccc}
+2 & +2 & +1 & +2 & +2\\
+2 & +1 & -1 & +2 & +2\\
-2 & +1 & +1 & +2 & +2\\
-2 & -1 & +2 & +2 & +2\\
+1 &\\
\vdots &\\
+1 &\\
-1 &\\
\end{array}\right]  
\end{displaymath} 

These four full-length words and the eight first letters of the rest of words from $V$, seven letters $+1$ and one $-1$, are the initial configurations of words for the computations corresponding to the case 1 of Corollary \ref{rozk} for $d=5$. Our task is to compute the missing cells of the matrix $M$ such that the rows in the resulting matrix form the polybox code $V$ with the following properties: $|V|=12$, the code $V$ does not contain a twin pair and, by Lemma \ref{ciecia}, for every $i\in [5]$ and every $l,s\in \{+1,-1,+2,-2\}, l\not\in \{s,s'\}$, the set $V^{i,l}\cup V^{i,s}$ contains an $i$-siblings. 

In the case 1, 2, 4, and 9 the initial configurations of words are the following (we write $M$ in the form of a table):    
\begin{displaymath}
\begin{tabular}{|l|r|}
\hline
\multicolumn{2}{|c|}{Case 1,\;$d=4,5$} \\ 
\hline
 $+2$ & $v^1$ \\
 $+2$ &  $v^2$ \\
 $-2$ & $v^3$ \\
 $-2$ & $v^4$  \\ \cline{1-2}
  $7\times +1$ & \; \\
  $1\times -1$ & \; \\
 \hline
 \end{tabular}\;\;
  \begin{tabular}{|l|r|}
 \hline 
 \multicolumn{2}{|c|}{Case 2,\;$d=4,5$} \\ 
\hline 
  $+2$ & $v^1$ \\
 $+2$ &  $v^2$ \\
 $-2$ & $v^3$ \\
 $-2$ & $v^4$  \\ \cline{1-2}
  $6\times +1$ & \; \\
  $2\times -1$ & \; \\
 \hline
 \end{tabular}\;\;
  \begin{tabular}{|l|r|}
 \hline  
  \multicolumn{2}{|c|}{Case 4,\;$d=4,5$} \\ 
\hline 
  $+2$ & $v^1$ \\
 $+2$ &  $v^2$ \\
 $-2$ & $v^3$ \\
 $-2$ & $v^4$  \\ \cline{1-2}
  $5\times +1$ & \; \\
  $3\times -1$ & \; \\
 \hline
 \end{tabular}\;\;
  \begin{tabular}{|l|r|}
 \hline 
  \multicolumn{2}{|c|}{Case 9,\;$d=4,5$} \\ 
\hline 
  $+2$ & $v^1$ \\
 $+2$ &  $v^2$ \\
 $-2$ & $v^3$ \\
 $-2$ & $v^4$  \\ \cline{1-2}
  $4\times +1$ & \; \\
  $4\times -1$ & \; \\
 \hline
 \end{tabular}
  \end{displaymath} 
where 
$$
v^1=+2+1+2(+2),\; v^2=+1-1+2(+2)
$$
$$
v^3=-1+1+2(+2),\; v^4=+1+2+2(+2),
$$
where $(+2)$ at the end of $v^i$ means that the letter $+2$ has to be placed at the fifth position in the words $v^i$, $i=1,2,3,4$, for the case $d=5$. This notation will be used also below.
  
Observe that in the cases 3, 5 and 10 of Corollary \ref{rozk}, by Statement \ref{stt} (d), it can be $|\{\breve{w}\cap \breve{v}\neq\emptyset: v\in V^{i,b}\}|=2$ and $|\{\breve{w}\cap \breve{v}\neq\emptyset: v\in V^{i,b'}\}|=3$ or $|\{\breve{w}\cap \breve{v}\neq\emptyset: v\in V^{i,b}\}|=2$ and $|\{\breve{w}\cap \breve{v}\neq\emptyset: v\in V^{i,b'}\}|=2$. Therefore, by Statement \ref{stt} (c),

\begin{displaymath}  
 \begin{tabular}{|l|r|}
\hline 
    \multicolumn{2}{|c|}{Case 3,\;$d=4,5$} \\ 
\hline  
  $+2$ & $w^1$ \\
 $+2$ &  $w^2$ \\
 $-2$ & $w^3$ \\
 $-2$ & $w^4$ \\
  $-2$ & $w^5$  \\ \cline{1-2}
    $6\times +1$ & \; \\
  $1\times -1$ & \; \\
  \hline
  \end{tabular}\;\;
\begin{tabular}{|l|r|}
  \hline 
     \multicolumn{2}{|c|}{Case 5,\;$d=4,5$} \\ 
\hline  
  $+2$ & $w^1$ \\
 $+2$ &  $w^2$ \\
 $-2$ & $w^3$ \\
 $-2$ & $w^4$ \\
  $-2$ & $w^5$  \\ \cline{1-2}
    $5\times +1$ & \; \\
  $2\times -1$ & \; \\
 \hline  
\end{tabular}\;\;
\begin{tabular}{|l|r|}
\hline
    \multicolumn{2}{|c|}{Case 10,\;$d=4,5$} \\ 
\hline 
  $+2$ & $w^1$ \\
 $+2$ &  $w^2$ \\
 $-2$ & $w^3$ \\
 $-2$ & $w^4$ \\
  $-2$ & $w^5$  \\ \cline{1-2}
    $4\times +1$ & \; \\
  $3\times -1$ & \; \\
 \hline
  \end{tabular}
\end{displaymath}
\begin{displaymath}  
 \begin{tabular}{|l|r|}
\hline 
    \multicolumn{2}{|c|}{Case 3,\;$d=4,5$} \\ 
\hline  
  $+2$ & $v^1$ \\
 $+2$ &  $v^2$ \\
 $-2$ & $v^3$ \\
   $-2$ & $v^4$  \\ \cline{1-2}
    $1\times -2$ & \; \\
    $6\times +1$ & \; \\
  $1\times -1$ & \; \\
  \hline
  \end{tabular}\;\;
\begin{tabular}{|l|r|}
  \hline 
     \multicolumn{2}{|c|}{Case 5,\;$d=4,5$} \\ 
\hline  
  $+2$ & $v^1$ \\
 $+2$ &  $v^2$ \\
 $-2$ & $v^3$ \\
  $-2$ & $v^4$  \\ \cline{1-2}
    $1\times -2$ & \; \\
    $5\times +1$ & \; \\
  $2\times -1$ & \; \\
 \hline  
\end{tabular}\;\;
\begin{tabular}{|l|r|}
\hline
    \multicolumn{2}{|c|}{Case 10,\;$d=4,5$} \\ 
\hline 
  $+2$ & $v^1$ \\
 $+2$ &  $v^2$ \\
 $-2$ & $v^3$ \\
  $-2$ & $v^4$  \\ \cline{1-2}
    $1\times -2$ & \; \\
    $4\times +1$ & \; \\
  $3\times -1$ & \; \\
 \hline
  \end{tabular}
\end{displaymath}

where
$$
w^1=+2+2+1(+2),\; w^2=+1+1-1(+2),
$$
$$
w^3=+1+1+2(+2),\;w^4=-1+1+1(+2),\; w^5=+2-1+1(+2)
$$
or
$$
w^1=+2+1+1(+2),\; w^2=+1+2-1(+2), 
$$
$$
w^3=+1+1+2(+2),\; w^4=-1+1+1(+2),\; w^5=+1-1-1(+2).
$$
\begin{uw} 
\label{wu}
{\rm In the codes above the sets of words  $\{w^1,w^2\}$ and  $\{w^3,w^4,w^5\}$ correspond to the codes $V$ and $W$ in Lemma \ref{24} for the case $|V|=2,|W|=3$.
Observe that the words $w^3,w^4,w^5$ can also take the form: $w^3=+1+1+2(+2), w^4=+1-1+1(+2), w^5=-1+2+1(+2)$, but permuting the letters at the first two positions we see that the polybox codes  $w^1=+2+2+1(+2), w^2=+1+1-1(+2), w^3=+1+1+2(+2), w^4=-1+1+1(+2), w^5=+2-1+1(+2)$ and  $w^1=+2+2+1(+2), w^2=+1+1-1(+2), w^3=+1+1+2(+2), w^4=+1-1+1(+2), w^5=-1+2+1(+2)$ are isomorphic. Clearly, we have to make the computations only for non-isomorphic codes. 

Moreover, the forms of the words $\{w^1,\ldots ,w^6\}$ are such as in Statement \ref{stt} (c) (see also Example \ref{p}). 

Similarly, in the case of the polybox codes $U=\{u^1,\ldots ,u^6\}$ and  $P=\{p^1,\ldots ,p^6\}$ which are given below, the codes $\{u^1,u^2\}$, $\{u^3,u^4,u^5,u^6\}$ and $\{p^1,p^2,p^3\}$, $\{p^4,p^5,p^6\}$ correspond to the codes $V$ and $W$ in Lemma \ref{24} for the cases $|V|=2,|W|=4$ and $|V|=3,|W|=3$, respectively. Also in these cases we  make the computations for non-isomorphic codes.}
\end{uw} 
 
By Statement \ref{stt} (c), (d) and the above remark, the cases 11 and 12 give the initial configurations of words:

\begin{displaymath}
\begin{tabular}{|l|r|}
 \hline 
  \multicolumn{2}{|c|}{Case 11,\;$d=4,5$} \\ 
\hline  
  $+2$ & $u^1$ \\
 $+2$ &  $u^2$ \\
 $-2$ & $u^3$ \\
 $-2$ & $u^4$ \\
 $-2$ & $u^5$ \\ 
  $-2$ & $u^6$  \\ \cline{1-2}
    $4\times +1$ & \; \\
  $2\times -1$ & \; \\
 \hline
 \end{tabular}\;\;
\begin{tabular}{|l|r|}
 \hline 
  \multicolumn{2}{|c|}{Case 11,\;$d=4,5$} \\ 
\hline  
  $+2$ & $w^1$ \\
 $+2$ &  $w^2$ \\
 $-2$ & $w^3$ \\
 $-2$ & $w^4$ \\ 
  $-2$ & $w^5$  \\ \cline{1-2}
    $1\times -2$ & \; \\
    $4\times +1$ & \; \\
  $2\times -1$ & \; \\
 \hline
 \end{tabular}\;\;
\begin{tabular}{|l|r|}
 \hline 
  \multicolumn{2}{|c|}{Case 11,\;$d=4,5$} \\ 
\hline  
  $+2$ & $v^1$ \\
 $+2$ &  $v^2$ \\
 $-2$ & $v^3$ \\ 
  $-2$ & $v^4$  \\ \cline{1-2}
    $2\times -2$ & \; \\
    $4\times +1$ & \; \\
  $2\times -1$ & \; \\
 \hline
 \end{tabular}
\end{displaymath}
\begin{displaymath}
\begin{tabular}{|l|r|}
\hline
   \multicolumn{2}{|c|}{Case 12,\;$d=4,5$} \\ 
\hline  
  $+2$ & $u^1$ \\
 $+2$ &  $u^2$ \\
 $-2$ & $u^3$ \\
 $-2$ & $u^4$ \\
 $-2$ & $u^5$ \\ 
  $-2$ & $u^6$  \\ \cline{1-2}
    $3\times +1$ & \; \\
  $3\times -1$ & \; \\
 \hline 
\end{tabular}\;\;
\begin{tabular}{|l|r|}
\hline
   \multicolumn{2}{|c|}{Case 12,\;$d=4,5$} \\ 
\hline  
  $+2$ & $w^1$ \\
 $+2$ &  $w^2$ \\
 $-2$ & $w^3$ \\
 $-2$ & $w^4$ \\ 
  $-2$ & $w^5$  \\ \cline{1-2}
    $1\times -2$ & \; \\
    $3\times +1$ & \; \\
  $3\times -1$ & \; \\
 \hline 
\end{tabular}\;\;
\begin{tabular}{|l|r|}
\hline
   \multicolumn{2}{|c|}{Case 12,\;$d=4,5$} \\ 
\hline  
  $+2$ & $v^1$ \\
 $+2$ &  $v^2$ \\
 $-2$ & $v^3$ \\ 
  $-2$ & $v^4$  \\ \cline{1-2}
    $2\times -2$ & \; \\
    $3\times +1$ & \; \\
  $3\times -1$ & \; \\
 \hline 
\end{tabular}
\end{displaymath}

where
$$
u^1=+2+2+2+1,\; u^2=+1+1+1-1,\; u^3=+1+1+1+2,
$$
$$
u^4=-1+1+1+1,\; u^5=+2-1+1+1,\; u^6=+2+2-1+1
$$
or
$$
u^1=+2+2+1+1,\; u^2=+1+1+2-1,\; u^3=+1+1+1+2
$$
$$
u^4=-1+1+1+1,\; u^5=+2-1+1+1,\; u^6=+1+1-1-1
$$
or
$$
u^1=+2+2+1(+2),\; u^2=+1+2-1(+2),\; u^3=+1+1+2(+2),
$$
$$
u^4=-1+1+1(+2),\; u^5=+2-1+1(+2),\; u^6=+1-1-1(+2).
$$

By Lemma \ref{33}, Statement \ref{stt} (c), (d) and Remark \ref{wu},

\begin{displaymath}
\begin{tabular}{|l|r|}
\hline
   \multicolumn{2}{|c|}{Case 13, $d=4$} \\ 
\hline  
  $+2$ & $p^1$ \\
 $+2$ &  $p^2$ \\
 $+2$ & $p^3$ \\
 $-2$ & $p^4$ \\
 $-2$ & $p^5$ \\ 
  $-2$ & $p^6$  \\ \cline{1-2}
    $3\times +1$ & \; \\
  $3\times -1$ & \; \\
 \hline 
\end{tabular}\;\;\;
\begin{tabular}{|l|r|}
\hline
   \multicolumn{2}{|c|}{Case 13, $d=4$} \\ 
\hline  
  $+2$ & $w^1$ \\
 $+2$ &  $w^2$ \\
 $+2$ & $w^3$ \\
 $-2$ & $w^4$ \\ 
  $-2$ & $w^5$  \\ \cline{1-2}
    $1\times -2$ & \; \\
    $3\times +1$ & \; \\
  $3\times -1$ & \; \\
 \hline 
\end{tabular}\;\;\;
\begin{tabular}{|l|r|}
\hline
   \multicolumn{2}{|c|}{Case 13, $d=4$} \\ 
\hline  
  $+2$ & $v^1$ \\
 $+2$ &  $v^2$ \\
 $-2$ & $v^3$ \\ 
  $-2$ & $v^4$  \\ \cline{1-2}
    $1\times +2$ & \; \\
    $1\times -2$ & \; \\
    $3\times +1$ & \; \\
  $3\times -1$ & \; \\
 \hline 
\end{tabular}
\end{displaymath}

where
$$
p^1=+2+1+1,\; p^2=+1+2-1,\; p^3=-1-1+2,
$$
\vspace{-6 mm}
$$
p^4=+2-1-1,\; p^5=-1+2+1,\; p^6=+1+1+2,
$$
or
$$
p^1=+1+2+2,\; p^2=-1-1+2,\; p^3=-1+1+1,
$$
\vspace{-6 mm}
$$
p^4=+2-1+2,\; p^5=+2+1+1,\; p^6=+1+1-1,
$$
or
$$
p^1=+1+2+2,\; p^2=-1-1+2,\; p^3=-1+1+1,
$$
\vspace{-6 mm}
$$
p^4=+2+2+1,\; p^5=+1+2-1,\; p^6=-1-1-1,
$$
or
$$
p^1=+1+2+2,\; p^2=-1-1+2,\; p^3=-1+1+1,
$$
\vspace{-6 mm}
$$
p^4=+2+2+1,\; p^5=+2-1-1,\; p^6=+1+1-1,
$$
or
$$
p^1=+1+1+2,\; p^2=+2-1-1,\; p^3=-1-1+1,
$$
\vspace{-6 mm}
$$
p^4=-1-1+2,\; p^5=+1+2-1,\; p^6=+1+1+1.
$$

\medskip
As we noted above, in the cases 6 and 7 we have $\breve{w}\cap \bigcup E(V^{i,a})\neq\emptyset$ and  $\breve{w}\cap \bigcup E(V^{i,a'})\neq\emptyset$. In this two cases we have $w_i=b$. Thus, by Lemma  \ref{rozk1}, 
$|\{\breve{w}\cap \breve{v}\neq\emptyset: v\in V^{i,a}\}|=5$ and $|\{\breve{w}\cap \breve{v}\neq\emptyset: v\in V^{i,a'}\}|=1$. Therefore, by Statement \ref{stt} (c),  

\begin{displaymath}
\begin{tabular}{|l|r|}
\hline
  \multicolumn{2}{|c|}{Case 6, $d=4,5$} \\ 
\hline 
  $+1$ & $q^1$ \\
 $+1$ &  $q^2$ \\
 $+1$ & $q^3$ \\
 $+1$ & $q^4$ \\
 $+1$ & $q^5$ \\ 
  $-1$ & $q^6$  \\ \cline{1-2}
    $3\times +2$ & \; \\
  $3\times -2$ & \; \\
\hline
\end{tabular}\;\;\;
\begin{tabular}{|l|r|}
\hline
  \multicolumn{2}{|c|}{Case 7, $d=4,5$} \\ 
\hline 
  $+1$ & $q^1$ \\
 $+1$ &  $q^2$ \\
 $+1$ & $q^3$ \\
 $+1$ & $q^4$ \\
 $+1$ & $q^5$ \\ 
  $-1$ & $q^6$  \\ \cline{1-2}
    $2\times +2$ & \; \\
  $4\times -2$ & \; \\
\hline
\end{tabular}
\end{displaymath}

where
\vspace{-4 mm}
$$
q^1=+1+1+1(+2),\; q^2=-1-1-1(+2),\; q^3=+2+1-1(+2),
$$
\vspace{-8 mm}
$$
q^4=-1+2+1(+2),\; q^5=+1-1+2(+2),\; q^6=+2+2+2(+2).
$$

\medskip
Finally, we consider the initial configurations of words corresponding to the case 8 of Corollary \ref{rozk}. In this case, by Lemma \ref{6}, besides computations for $d=4,5$ we have to make  computations for $d=6$. 

Observe first that, by Statement \ref{stt} (d),  $V^{i,b}\sqsubseteq W^{i,b}$ and $V^{i,b'}\sqsubseteq W^{i,b'}$, and thus, by Statement \ref{stt} (a),  $|W^{i,b}|\geq 5$ and $|W^{i,b'}|\geq 5$, and since, by Lemma \ref{ciecia}, $|W^{i,a}|\geq 1$ and $|W^{i,a'}|\geq 1$, we have $|W^{i,b}|=5$, $|W^{i,b'}|=5$ and $|W^{i,a}|=1$, $|W^{i,a'}|=1$. Clearly, $W^{i,a}\sqsubseteq V^{i,a}$ and $W^{i,a'}\sqsubseteq V^{i,a'}$, and thus the structure of  $V^{i,a}$ and  $V^{i,a'}$ are such as in Statement \ref{stt} (a).


We do not know what is a relation between the polybox codes $V^{i,a}$ and  $V^{i,a'}$. Therefore, we fix the structure of $V^{i,a}$, say that  $(V^{i,a})_i=\{q^1,\ldots, q^5\}$ for $d=4,5$ and $(V^{i,a})_i=\{q^1+2,\ldots, q^5+2\}$ for $d=6$, and $V^{i,a'}$ will take all possible forms, i.e., $(V^{i,a'})_i=\{x^1_\sigma,\ldots, x^5_\sigma\}$, where
$$
x^1=l_1l_2l_3(u_4)(u_5),\; x^2=l_1'l_2'l_3'(u_4)(u_5),\; x^3=s_1l_2l_3'(u_4)(u_5)
$$
$$
x^4=l_1's_2l_3(u_4)(u_5),\; x^5=l_1l_2's_3(u_4)(u_5),
$$
$s_1,s_2,s_3,l_1,l_2,l_3,u_4,u_5$ range over the set $\{+1,-1,+2,-2\}$, $s_i\not\in \{l_i,l_i'\}$ for $i=1,2,3$, $\sigma$ is a permutation of the set $[d]$ and $v_\sigma=v_{\sigma(1)}\cdots v_{\sigma(d)}$, where $v=v_1\cdots v_d$. Now, for every selection of  letters $s_1,s_2,s_3,l_1,l_2,l_3,u_4,u_5\in \{+1,-1,+2,-2\}$ and a permutation $\sigma$ we have the initial configuration 

\medskip
\begin{displaymath}
\begin{tabular}{|l|r|}
\hline
  \multicolumn{2}{|c|}{Case 8,\;$d=4,5,6$} \\ 
\hline 
  $+1$ & $q^1$ \\
 $+1$ &  $q^2$ \\
 $+1$ & $q^3$ \\
 $+1$ & $q^4$ \\
 $+1$ & $q^5$\\
 $-1$ & $x^1_\sigma$ \\
 $-1$ &  $x^2_\sigma$ \\
 $-1$ & $x^3_\sigma$ \\
 $-1$ & $x^4_\sigma$ \\ 
   $-1$ & $x^5_\sigma$  \\ \cline{1-2}
    $1\times +2$ & \; \\
  $1\times -2$ & \; \\
\hline
\end{tabular}
\end{displaymath}

\medskip
The  case 8 generates quite a large  package of the initial configurations of words, but each of them contains as many as ten words, which makes every single computation very quick.  

Below we present an algorithm using in the computations.

Let $S=\{+1,-1,+2,-2\}$, and let $S^d_{+1}$, $S^d_{-1}$ and $S^d_{-2}$ denote the sets of all words of length $d$, for $d=4,5,6$, having the letter $+1$, $-1$ or $-2$ at the first position, respectively. By $M^d_r(n_{+1},n_{-1},n_{+2})$ we denote the set of $r$ full-length words in the initial configuration having $n_{+1},n_{-1}$ and $n_{+2}$ letters $+1,-1$ and $+2$ at the first position, respectively.

\textbf{Algorithm} 

\noindent
\textbf{Input:} The set $M^d_r(n_{+1},n_{-1},n_{+2})$

\noindent
\textbf{Output:} A polybox code $V$ having $12$ words without twin pairs such that for every $i\in[d]$ and every $l,s\in\{+1,-1,+2,-2\}$, $l\notin\{s,s'\}$, the set $V^{i,l}\cup V^{i,s}$ contains an $i$-siblings

\medskip

\begin{flushleft}
\begin{tabular}{r p{12cm}}
1:& $V:=\emptyset$\\
2:& $W:=M^d_r(n_{+1},n_{-1},n_{+2})$\\
3:& $k:=r$\\
4:& \textbf{repeat}\\
5:&  \hspace{0.5cm}\textbf{while} $k<12$ \textbf{do}\\
6:&    \hspace{1cm}select a word $v$ from the set $S^d_{+1}\cup S^d_{-1}\cup S^d_{-2}$ which is dichotomous to every word of $W$ and does not form a twin pair with any word of $W$\\
7:&    \hspace{1cm}$W:=W\cup\{v\}$\\
8:&    \hspace{1cm}$k:=k+1$\\
9:&  \hspace{0.5cm}\textbf{end while}\\
10:&  \hspace{0.5cm}\textbf{if} $|W^{1,+1}|=n_{+1},|W^{1,-1}|=n_{-1}$, $|W^{1,-2}|=n_{-2}$ and for every $i\in[d]$ and  $l,s\in\{+1,-1,+2,-2\}$, $l\notin\{s,s'\}$, the set $W^{i,l}\cup W^{i,s}$ contains an $i$-siblings \textbf{then}\\
11:&   \hspace{1cm}$V:=W$\\
12:&   \hspace{1cm}\textbf{return} V\\
13:&   \hspace{0.5cm}\textbf{end if}\\
14:&   \hspace{0.5cm}$W:=M^d_r(n_{+1},n_{-1},n_{+2})$\\
15:&   \hspace{0.5cm}$k:=r$\\
16:&  \textbf{until} for every  $m\in \{+1,-1,-2\}$ all possible subsets of $S^d_m$ with $n_m$ elements are selected
\end{tabular}
\end{flushleft}

\bigskip
For every two words $v$ and $w$ by $vw$ we denote a concatenation of $v$ and $w$. If $A$ and $B$ are sets of words, then
$$
AB=\{vw:v\in A,\,\, w\in B\}.
$$

The results of the computations are given in the following theorem: 

\begin{tw}
\label{s12}
Let $l,s\in \{a,a',b,b'\}, l\not\in \{s,s'\}$, and let  
$$
A_1=\{ls'\},\, B_1=\{l'l'\},\, C_1=\{ls, l'l\},\, A_1^c=\{ss, l's', s's\},\, B_1^c=\{ll', sl, s'l\},
$$
$$
A_2=\{s'l\},\, B_2=\{ss\},\, C_2=\{ss', s'l'\},\, A_2^c=\{ll', l'l', sl\},\, B_2^c=\{ls', l's', s's\},
$$
$$
W_1=C_1C_1\cup A_1A_1^c\cup B_1B_1^c\cup A_1^cB_1\cup B_1^cA_1,
$$
$$
W_2=C_2C_2\cup A_2^cA_2\cup B_2^cB_2\cup A_2B_2^c\cup B_2A_2^c.
$$
If $V,W\subset S^d$, where $d\in \{4,5,6\}$, are disjoint and equivalent polybox codes without twin pairs having twelve words and  $V$ is extensible to a partition code, then there is a set $A=\{i_1<i_2<i_3<i_4\}\subset [d]$ such that 
$$
(V)_{A^c}=W_1\setminus W_1\cap W_2,\;\;\;   (W)_{A^c}=W_2\setminus W_1\cap W_2,
$$
and  $(V)_{A}=(W)_{A}=\{(p)_A\}$ for some $p\in S^d$. The representation of $V$ and $W$ is given up to an isomorphism.
\end{tw}   
\proof 
To show that $V$ and $W$ do not contain a twin pair observe that the partition code $W_1$ contains four twin pairs
$$
ls'ss,\,\, ls's's;\;\; l'l's'l,\,\, l'l'sl;\;\; ssl'l',\,\, s'sl'l';\;\; slls',\,\, s'lls',
$$
and the partition code $W_2$ also contains four twins  
$$
\mathbf{ ls'ss},\,\, l's'ss;\;\; \mathbf{l'l's'l},\,\, ll's'l;\;\;
\mathbf{ ssl'l'},\,\, ssll';\;\; \mathbf{ s'lls'},\,\, s'll's',
$$
where the marked words form the set $W_1\cap W_2$. Therefore, the polybox codes $V=W_1\setminus W_1\cap W_2$ and $W=W_2\setminus W_1\cap W_2$ are disjoint and equivalent, they do not contain twin pairs  and $|V|=|W|=12$.

By Corollary \ref{rozk}, the distributions of the letters in $V$ cannot be other than the distributions $1-13$ of this corollary. Thus, the code $V$ has to contain at least one initial configuration from the list of the initial configurations of words given in the tables in Section 5.1.

The computations showed that the initial configurations corresponding to the cases $1-7$ and $9-12$ cannot be completed to a twin pair free polybox code $V$ with twelve words such that for every $i\in[d]$ and $l,s\in\{a,a',b,b'\}$, $l\notin\{s,s'\}$, the set $V^{i,l}\cup V^{i,s}$ contains an $i$-siblings (compare Lemma \ref{ciecia}).
Similarly, this cannot be done in the case $8$ for $d=5$ and 6. The only configurations that can be completed to a polybox code $V$ with the above properties are the case 8 for $d=4$ and the configurations given in the last two tables for the case 13. The structures of all these complemented codes $V$ (for $d=4$ we have $(V)_{A^c}=V$ as $A^c=\emptyset$) with 12 words are, up to an isomorphism, as stated in the theorem. The polybox code $W$ is a complementation  of the words $W_1\cap W_2=\{ls'ss, l'l's'l, ssl'l', slls'\}$ to a partition codes. It can be easily check that there is exactly one such complementation which is disjoint with $V$, and its structure is such as given in the theorem. (The words $ls'ss, l'l's'l, ssl'l', slls'$ can be complemented to a partition code exactly in two ways: $V$ and $W$). 

Since for $d=5$ and 6 all initial configurations cannot be  completed to a twin pair free polybox code with twelve words such that for every $i\in[d]$ and $l,s\in\{a,a',b,b'\}$, $l\notin\{s,s'\}$, the set $V^{i,l}\cup V^{i,s}$ contains an $i$-siblings, it follows that for the above two dimensions we have, by Lemma \ref{ciecia}, $(V)_{A}=(W)_{A}=\{(p)_A\}$ for some $p\in S^d$, and $(V)_{A^c}$, $(W)_{A^c}$ are such as in dimension four. 
\hfill{$\square$}

\medskip
\begin{uw}
{\rm The representation of $V$ and $W$ given in Theorem \ref{s12} can be found in \cite{KP1}. The codes $V$ and $W$ were used by Lagarias and Shor\cite{LS1,LS2} and later on by Mackey\cite{M} to construct the counterexamples to Keller's cube tiling conjecture. In the context of this conjecture one of these codes was given first by Corr\'adi and Szab\'o  in \cite{CS2}, as an example of the maximal clique in a 4-dimensional Keller graph.  
}
\end{uw}

\section{Twin pairs in cube tilings of $\er^7$}
  
From Theorem \ref{12} and \ref{s12} we obtain the following 

\begin{tw}
\label{112}
Let $U\subset S^7$ be a partition code. If there are $i\in [7]$ and $l\in S$ such that $|U^{i,l}|\leq 12$, then there is a twin pair in $U$.
\end{tw}
\proof Since $\bigcup E(U^{i,l}\cup U^{i,l'})$ is an $i$-cylinder in $(ES)^7$, the codes $U^{i,l}$ and $U^{i,l'}$ are equivalent. If  $U^{i,l}$ or $U^{i,l'}$ contains a twin pair, then clearly $U$ does. Thus, we assume that these two codes do not contain a twin pair. If $(U^{i,l})_i\cap (U^{i,l'})_i\neq\emptyset$ and  $(v)_i\in (U^{i,l})_i\cap (U^{i,l'})_i$, then the words $w\in U^{i,l}$ and $p\in U^{i,l'}$ such that $(v)_i=(w)_i=(p)_i$ are a twin pair. Therefore we may assume that  $(U^{i,l})_i$ and $(U^{i,l'})_i$ are disjoint and do not contain a twin pair. It follows from Theorem \ref{12} that $|U^{i,l}|=12$, and Theorem \ref{s12} precisely describes the structure of the codes $(U^{i,l})_i$ and $(U^{i,l'})_i$ (clearly, by (\textbf{S}) in Section 3.1, the codes $(U^{i,l})_i$ and $(U^{i,l'})_i$ are extensible to partition codes). In Theorem \ref{s12} we take $A=\{1,2,3,4\}, (p)_A=ll$ and $(U^{i,l})_i=V$. Without loss of generality we can take $i=7$. We consider only the code $U^{7,l}$ which has the form:

\begin{displaymath}
\begin{array}{rrrrrrr}
l&s&l&s&l&l&l\\
l&s&l'&l&l&l&l\\
l'&l&l&s&l&l&l\\
l'&l&l'&l&l&l&l\\
l&s'&l'&s'&l&l&l\\
l&s'&s'&s&l&l&l\\
l'&l'&l&l'&l&l&l\\
l'&l'&s&l&l&l&l\\
l'&s'&l'&l'&l&l&l\\
s'&s&l'&l'&l&l&l\\
l&l'&l&s'&l&l&l\\
s&l&l&s'&l&l&l\\
\end{array}
\end{displaymath}

We choose four words from  $U^{7,l}$:
\begin{displaymath}
\begin{array}{rrrrrrr}
v=l&s&l&s&l&l&l\\
u=l&s&l'&l&l&l&l\\
p=l&s'&l'&s'&l&l&l\\
q=s'&s&l'&l'&l&l&l\\
\end{array}
\end{displaymath}
Let $W\subset U^{4,l'}$ be the set of all $w\in U^{4,l'}$ such that $(u)_4\sqsubseteq (W)_4$ and $(\breve{w})_4\cap  (\breve{u})_4\neq\emptyset$. Since $U$ is a partition code, $W\neq\emptyset$. Clearly, for every $i\in [7], i\neq 4,$ and $w\in W$ we have $w_i\neq u_i'$, for otherwise $(\breve{w})_4\cap  (\breve{u})_4=\emptyset$ which contradicts the definition of $W$. Every $w\in W$ is dichotomous to the words $v,q$ and $p$, and therefore $w_1=w_2=s$ and  $w_3=w_4=l'$ for every $w\in W$. Then
$$
lll\sqsubseteq (W)_{\{1,2,3,4\}}.
$$   
If $(W)_{\{1,2,3,4\}}$ does not contain a twin pair, then its structure is such as in Statement \ref{stt} (a). Without loss of generality we may assume that 
$$
(W)_{\{1,2,3,4\}}=\{s_1s_2s_3,s_1's_2's_3',ls_2's_3,s_1ls_3',s_1's_2l\},
$$
where $s_i\not\in \{l,l'\}$ for $i=1,2,3$. But then $ssl'l's_1's_2l\in U^{7,l}$ which is not true. Therefore, the code $(W)_{\{1,2,3,4\}}$ contains a twin pair, and hence $W$ contains a twin pair.
\hfill$\square$

\begin{wn}
Let $[0,1)^7+T$ be a cube tiling of $\er^7$. If $|L(T,x,i)|=5$ for some $x\in \er ^7$ and $i\in [7]$, then there is a twin pair in $[0,1)^7+T$.
\end{wn} 
\proof  As it was noted at the beginning of  Section 2, the set of boxes $\ka F_x=\{([0,1)^7+t)\cap ([0,1]^7+x)\neq\emptyset:t\in T\}$ is a minimal partition of $[0,1]^7+x$. 
Let $U$ be a partition code such that $\ka F_x$ is an exact realization of $U$ (compare \cite[Theorem 4.2]{Kis}). Since  $|L(T,x,i)|=5$, we have $U=U^{i,l_1}\cup U^{i,l_1'}\cup \cdots \cup U^{i,l_5}\cup U^{i,l_5'}$ and   $|U^{i,l_j}|\leq 12$ for some $j\in [5]$. By Theorem \ref{112}, there is a twin pair in $U$, and consequently there is a twin pair in $\ka F_x$. Then the tiling  $[0,1)^7+T$ contains a twin pair.
\hfill{$\square$}

\medskip
From the result of Debroni et al., \cite[Corollary 4.2]{Kis} and the above corollary we obtain the following

\begin{wn}
If there is a counterexamples to Keller's conjecture in dimension seven, then $|L(T,x,i)|\in \{3,4\}$ for some $x\in \er ^7$ and $i\in [7]$.
\end{wn} 

\medskip
In \cite{Kis} we extended the notion of a $d$-dimensional Keller graph:  if $S$ is an alphabet with a complementation, then {\it a d-dimensional Keller graph on the set $S^d$} is the graph in which two vertices $u,v\in S^d$ are adjacent if they are dichotomous but do not form a twin pair. 

From Theorem \ref{112} we obtain the following 
\begin{wn}
Every clique in a $7$-dimensional Keller graph on $S^7$ which contains at least five vertices $u^1,\ldots ,u^5$ such that $u^n_i\not\in \{u^m_i,(u^m_i)'\}$ for some $i\in [7]$  and every $n,m\in \{1,...,5\},n\neq m$, has less than $2^7$ elements. 
\end{wn}
\proof Assume on the contrary that there is a clique $U$ containing vertices $u^1,\ldots ,u^5$ and $|U|=2^7$. Thus, $U$ is a partition code without twin pairs.  Since  $u^n_i\not\in \{u^m_i,(u^m_i)'\}$ for every $n,m\in \{1,...,5\},n\neq m$, 
it follows that $|U^{i,u^m_i}|\leq 12$ for some $m\in [5]$. By Theorem \ref{112}, there is a twin pair in $U$, a contradiction.  
\hfill{$\square$}

\begin{uwi} {\rm What next? Our computer experiments made together with our colleague Krzysztof Przes\l awski show that there is a chance that the case $|L(T,x,i)|=4$ can be resolve using the same methods as the case $|L(T,x,i)|=5$, but a scale of the computations will be much more larger than that for $|L(T,x,i)|=5$. Moreover, to obtain initial configurations for the computations we need some new results on the rigidity of polyboxes.  

}
\end{uwi}

\end{document}